\documentclass[twoside]{article}
\usepackage{graphicx, epstopdf, verbatim, fancyhdr, indent, wrapfig}

\usepackage{amsmath}
\input amssym.def
\input amssym.tex
\usepackage{bbm}

\skip\footins 22pt plus 2pt minus 2pt
\newdimen\jot
\def\mqth{\mathsurround=0pt }
\dimendef\dimenq=0
\def\openup{\afterassignment\qpenup\dimenq=}
\def\qpenup{\advance\lineskip\dimenq
  \advance\baselineskip\dimenq \advance\lineskiplimit\dimenq}
\def\eqalign#1{\,\vcenter{\openup1\jot \mqth
  \ialign{\strut\hfil$\displaystyle{##}$&$\displaystyle{{}##}$\hfil
  \crcr#1\crcr}}\,}
\newif\ifdtqp
\def\displqy{\global\dtqptrue \openup1\jot \mqth
  \everycr{\noalign{\ifdtqp \global\dtqpfalse
     \vskip-\lineskiplimit \vskip\normallineskiplimit
     \else \penalty\interdisplaylinepenalty \fi}}}
\def\displaylines#1{\displqy
  \halign{\hbox to\displaywidth{$\hfil\displaystyle##\hfil$}\crcr
  #1\crcr}}
\newskip\centerinq \centerinq=0pt plus 1000pt minus 1000pt

\newcommand{\eps}{\varepsilon}
\newcommand{\ph}{\varphi}
\newcommand{\dast}{{\displaystyle{\ast}}}
\newcommand{\reals}{{\mathbbm{R}}}

\newcommand{\zero}{{\bf 0}}

\newcommand{\csch}{\mathop{\rm csch}\nolimits}

\newcommand{\Circ}{\mathop{\rm Circ}\nolimits}
\newcommand{\Lk}{\mathop{\rm Lk}\nolimits}
\newcommand{\Tw}{\mathop{\rm Tw}\nolimits}
\newcommand{\Wr}{\mathop{\rm Wr}\nolimits}
\newcommand{\BS}{\mathop{\rm BS}\nolimits}
\newcommand{\VF}{\mathop{\rm VF}\nolimits}
\newcommand{\He}{\mathop{\raise.02ex\hbox{\rm H}}\nolimits}
\newcommand{\K}{\mathop{\raise.02ex\hbox{\rm K}}\nolimits}
\newcommand{\G}{\mathop{\raise.02ex\hbox{\rm G}}\nolimits}
\newcommand{\CK}{\mathop{\rm CK}\nolimits}
\newcommand{\AK}{\mathop{\rm AK}\nolimits}
\newcommand{\HK}{\mathop{\rm HK}\nolimits}
\newcommand{\vol}{\mathop{\rm vol}\nolimits}
\newcommand{\area}{\mathop{\rm area}\nolimits}
\newcommand{\curl}{\mathop{\rm curl}\nolimits}
\newcommand{\x}{{\bf x}}
\newcommand{\y}{{\bf y}}
\newcommand{\ye}{{\bf y}_{\eps}}
\newcommand{\ys}{{\bf y}_{\sigma}}
\newcommand{\as}{{\alpha_\sigma}}
\newcommand{\aes}{{\alpha_\eps(s)}}
\newcommand{\vv}{{\bf v}}
\newcommand{\w}{{\bf w}}

\newcommand{\T}{{\bf T}}
\newcommand{\N}{{\bf N}}
\newcommand{\B}{{\bf B}}

\newcommand{\A}{{\bf A}}
\renewcommand{\i}{{\bf i}}
\renewcommand{\j}{{\bf j}}

\renewcommand{\u}{{\bf u}}
\newcommand{\n}{{\bf n}}

\newcommand{\xh}{{\hat{\x}}}

\def\ip#1#2{\langle#1,#2\rangle}

\begin{document}
\thispagestyle{empty}

\large
\LARGE
\centerline{\bf Linking, twisting, writhing and helicity}
\centerline{\bf on the 3-sphere and in hyperbolic 3-space}
\large

\bigskip

\centerline{\sc Dennis DeTurck and Herman Gluck}

\addtolength{\baselineskip}{1pt}

In the first paper of this series, ``\it Electrodynamics and the Gauss Linking Integral on the
$3$-sphere and in hyperbolic $3$-space\rm,'' we developed a steady-state version of classical
electrodynamics in these two spaces, including explicit formulas for the vector-valued
Green's operator, explicit formulas of Biot-Savart type for the magnetic field, and a
corresponding Amp\`ere's Law contained in Maxwell's equations, and then used these to
obtain explicit integral formulas for the linking number of two disjoint closed curves.

In this second paper, we obtain integral formulas for twisting, writhing and helicity, and prove
the theorem \sc link = twist + writhe \rm on the 3-sphere and in hyperbolic 3-space.
We then use these results to derive upper bounds for the helicity of vector fields and lower
bounds for the first eigenvalue of the curl operator on subdomains of these two spaces.

An announcement of these results, and a hint of their proofs, can be found in the Math
ArXiv, math.GT/0406276, while an expanded version of the first paper, with full proofs, can
be found at math.GT/0510388. 

The flow of this paper is indicated by the following list of sections. The first two are devoted
to a summary of information from the preceding paper.
\begin{enumerate}
\item  Linking integrals in $\reals^3$, $S^3$ and $H^3$.
\item Magnetic fields in $\reals^3$, $S^3$ and $H^3$.
\item Link, twist and writhe in $S^3$ and $H^3$.
\item Proof scheme for \sc link = twist + writhe\rm.
\item Some geometric formulas on $S^3$.
\item Proof of \sc link = twist + writhe \rm in $S^3$.
\item Proof of \sc link = twist + writhe \rm in $H^3$.
\item Helicity of vector fields on $S^3$ and $H^3$.
\item Upper bounds for helicity in $\reals^3$, $S^3$ and $H^3$.
\item Hodge decomposition of vector fields.
\item Spectral geometry of the curl operator in $\reals^3$, $S^3$ and $H^3$.
\end{enumerate}

The integral formulas in this paper contain vectors lying in different tangent spaces; in
non-Euclidean settings these vectors must be moved to a common location to be
combined.

In $S^3$ regarded as the group of unit quaternions, equivalently as $SU(2)$, the
differential $L_{\y\x^{-1}}$ of left translation by $\y\x^{-1}$ moves tangent vectors from $\x$
to $\y$. In either $S^3$ or $H^3$, parallel transport $P_{\y\x}$ along the geodesic segment
from $\x$ to $\y$ also does this. As a result, we get three versions for each of the formulas
that appear in the theorems below.

\bigskip

\noindent\bf 1.\ Linking integrals in $\reals^3$, $S^3$ and $H^3$.\rm

Let $K_1 = \{\x(s)\}$ and $K_2=\{\y(t)\}$ be disjoint oriented smooth closed curves in
either Euclidean 3-space $\reals^3$, the unit 3-sphere $S^3$, or hyperbolic 3-space
$H^3$, and let $\alpha(\x,\y)$ denote the distance from $\x$ to $\y$.

%%% FIGURE 1 %%%
\begin{figure}[h!]
\center{\includegraphics[height=160pt]{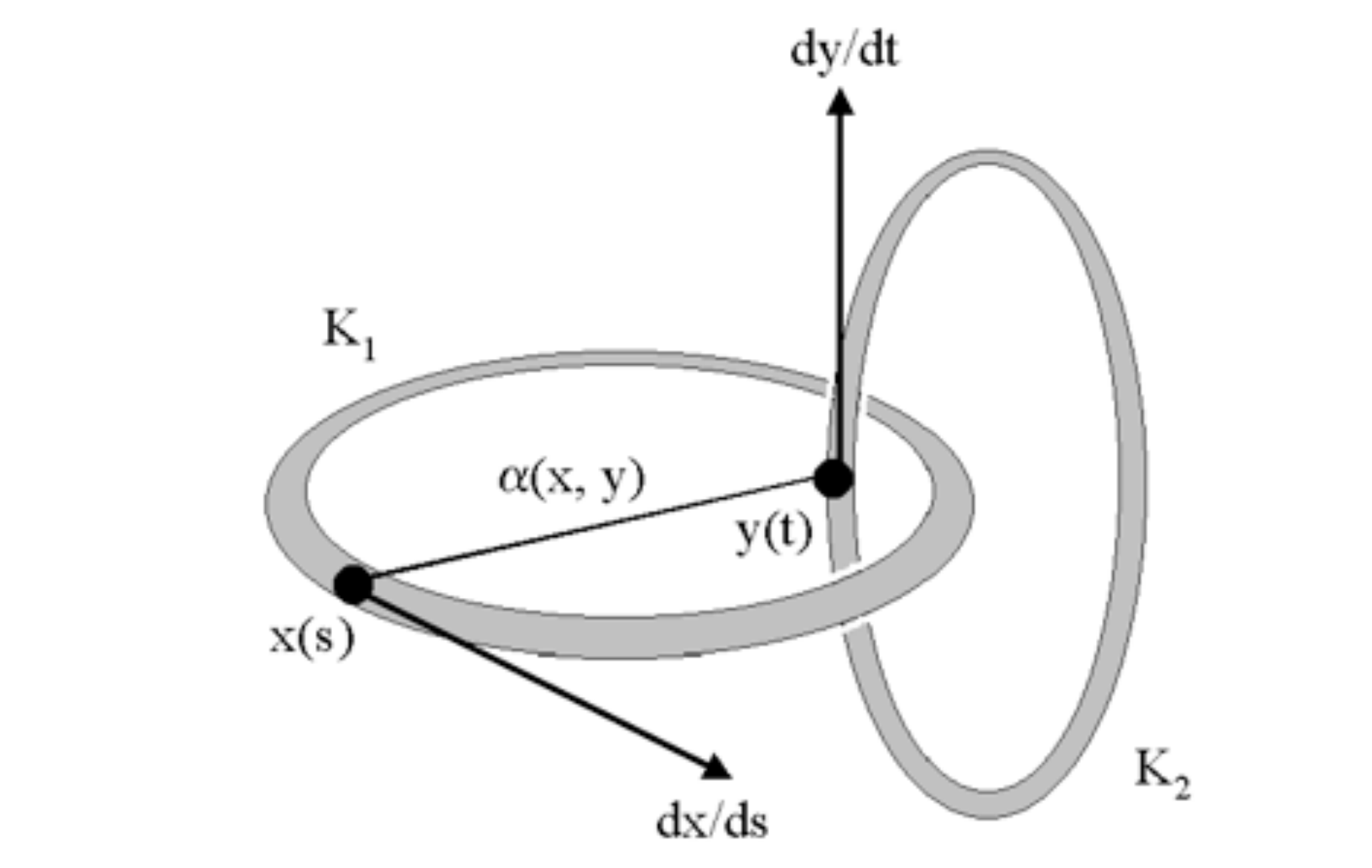}
\caption{Two linked curves}}
\end{figure}
%%%%%%%%%%%%%%%%%%%%

Carl Friedrich Gauss, in a half-page paper dated January 22, 1833, gave an integral
formula for the linking number in Euclidean 3-space,
$$\Lk(K_1,K_2) = \int_{K_1\times K_2} \frac{d\x}{ds}\times
\frac{d\y}{dt}\cdot\frac {\x-\y}{ 4\pi|\x-\y|^3}\,ds\,dt.$$
It will be convenient for us to write this as
$$\Lk(K_1,K_2) = \int_{K_1\times K_2} \frac{d\x}{ds}\times
\frac{d\y}{dt}\cdot \nabla_\y\ph(\x,\y)\,ds\,dt,$$
where $\ph(\alpha)=1/(4\pi\alpha)$, and where we use $\ph(\x,\y)$ as an abbreviation for
$\ph(\alpha(\x,\y))$. The subscript $\y$ in the expression $\nabla_\y\ph(\x,\y)$ tells
us that the differentiation is with respect to the $\y$ variable.

The following theorem from our first paper gives the corresponding linking integrals
on the 3-sphere and in hyperbolic 3-space. Since the location of the tangent vectors
is now important, we note that the vector $\nabla_\y\ph(\x,\y)$ is located at the point $\y$.

\medskip

\noindent\bf Theorem 1.1\rm.\ \sc Linking integrals in $S^3$ and $H^3$.

\rm(1) \it On $S^3$ in left-translation format:
$$\Lk(K_1,K_2) = \int_{K_1\times K_2} L_{\y\x^{-1}}\frac{d\x}{ds}\times
\frac{d\y}{dt}\cdot \nabla_\y\ph(\x,\y)\,ds\,dt
-\frac{1}{4\pi^2}\int_{K_1\times K_2} L_{\y\x^{-1}}\frac{d\x}{ds}
\cdot\frac{d\y}{dt}\,ds\,dt,$$
\indent where $\ph(\alpha)= (\pi-\alpha)\cot\alpha/(4\pi^2)$.

\rm(2) \it On $S^3$ in parallel transport format:
$$\Lk(K_1,K_2) = \int_{K_1\times K_2} P_{\y\x}\frac{d\x}{ds}\times
\frac{d\y}{dt}\cdot \nabla_\y\ph(\x,\y)\,ds\,dt,$$
\indent where $\ph(\alpha)=(\pi-\alpha)\csc\alpha/(4\pi^2)$.

\rm(3) \it On $H^3$ in parallel transport format:
$$\Lk(K_1,K_2) = \int_{K_1\times K_2} P_{\y\x}\frac{d\x}{ds}\times
\frac{d\y}{dt}\cdot \nabla_\y\ph(\x,\y)\,ds\,dt,$$
\indent where $\ph(\alpha)=\csch\alpha/(4\pi)$.\rm

\medskip

Greg Kuperberg (2008) obtained, independently and by a totally different
argument, an expression equivalent to formula (2) above.

The kernel functions used here have the following significance.

In Gauss's linking integral, the function $-\ph(\alpha)=-1/(4\pi\alpha)$, where $\alpha$
is distance from a fixed point, is the fundamental solution of the Laplacian in $\reals^3$,
$$-\Delta \ph=\delta.$$
Here $\delta$ is the Dirac $\delta$-function.

In formula (1), the function $-\ph(\alpha)= -(\pi-\alpha)\cot\alpha/(4\pi^2)$, is the
fundamental solution of the Laplacian on $S^3$,
$$-\Delta\ph = \delta - \frac{1}{2\pi^2}.$$
Since the volume of $S^3$ is $2\pi^2$, the right-hand side has average value zero.

In formula (2), the function $-\ph(\alpha) = -(\pi-\alpha)\csc\alpha/(4\pi^2)$ is the
fundamental solution of a shifted Laplacian on $S^3$,
$$-(\Delta\ph-\ph)=\delta.$$

In formula (3), the function $-\ph(\alpha)= -\csch\alpha/(4\pi)$ is the fundamental
solution of a shifted Laplacian on $H^3$,
$$-(\Delta\ph+\ph)=\delta.$$

Our proof of the formula \sc link = twist + writhe \rm will depend on the asymptotic
properties of $\ph$ at its singularity. For example, in the case of $S^3$ in
parallel transport format,
$$\ph(\alpha)=\frac{1}{4\pi^2}(\pi-\alpha)\csc(\alpha)= \frac{1}{4\pi\alpha}-\frac{1}{4\pi^2}+
\frac{1}{24\pi}\alpha -
\frac{1}{24\pi^2}\alpha^2 + \alpha^3
f(\alpha),$$
where $f(\alpha)$ is bounded and smooth. Likewise
$$\ph'(\alpha)=-\frac{1}{4\pi\alpha^2}+\frac{1}{24\pi}-
\frac{1}{12\pi^2}\alpha + \frac{7}{480\pi}\alpha^2
-\alpha^3g(\alpha)$$
and
$$\ph''(\alpha)=\frac{1}{2\pi\alpha^3}-\frac{1}{12\pi^2}+\frac{7}{240\pi}\alpha
-\frac{7}{120\pi^2}\alpha^2+\alpha^3h(\alpha),$$
where $g$ and $h$ are also bounded and smooth. Note that
$\ph$ has no singularity at $\alpha=\pi$,
in fact, $\ph$ is smooth and even around $\alpha=\pi$:
$$\ph(\alpha)=\frac{1}{4\pi^2}+\frac{1}{24\pi^2}(\alpha-\pi)^2+\frac{7}{1440\pi^2}(\alpha-\pi)^4+\cdots$$
near $\pi$. This implies that $\nabla_\y\ph(\alpha(\x,\y))$ exists and is
zero when $\y$ is the antipodal point of $x$, even though
$\nabla_\y\alpha$
is not defined there. Because of this, the functions $f(\alpha(\x,\y))$,
$g(\alpha(\x,\y))$ and $h(\alpha(\x,\y))$ defined above
are defined, smooth and bounded for \it all \rm $\x$ and $\y$ such that
$\x\ne \y$.

Because we do not need so many terms of these expansions, we will simply write:
$$\ph(\alpha)=\frac{1}{4\pi\alpha}+f(\alpha), \qquad \ph'(\alpha)=-\frac{1}{4\pi\alpha^2} +
g(\alpha), \qquad \ph''(\alpha)=\frac{1}{2\pi\alpha^3}+h(\alpha),$$
where these new functions $f$, $g$ and $h$ are bounded and smooth everywhere on $S^3$.

Similar calculations show that, for $\ph(\alpha)= \csch\alpha/(4\pi)$ on $H^3$,
we again have
$$\ph(\alpha)=\frac{1}{4\pi\alpha}+f(\alpha), \qquad \ph'(\alpha)=-\frac{1}{4\pi\alpha^2} +
g(\alpha), \qquad \ph''(\alpha)=\frac{1}{2\pi\alpha^3}+h(\alpha),$$
where these latest functions $f$, $g$ and $h$ are bounded and smooth everywhere on $H^3$.

\vfill
\eject

\noindent\bf 2.\ Magnetic fields in $\reals^3$, $S^3$ and $H^3$.\rm

In Euclidean 3-space $\reals^3$, the classical convolution formula of Biot and Savart
gives the magnetic field $\BS(\vv)$ of a compactly supported current flow $\vv$:
$$\BS(\vv)(\y)=\int_{\reals^3} \vv(\x)\times \frac{\y-\x}{4\pi|\y-\x|^3}\,d\x.$$
For simplicity, we write $d\x$ to mean $d\vol_\x$.

The Biot-Savart formula can also be written as
$$\BS(\vv)(\y)=\int_{\reals^3} \vv(\x)\times\nabla_\y\ph_0(\x,\y)\,d\x,$$
where $\ph_0(\alpha)=-1/(4\pi\alpha)$ is the fundamental solution of the Laplacian
in $\reals^3$.

In $\reals^3$, if we start with a smooth, compactly supported current flow $\vv$, then its
magnetic field $\BS(\vv)$ is a smooth vector field (although not in general compactly
supported) which has the following properties:

\begin{enumerate}
\item[(1)] It is divergence-free, $\nabla\cdot\BS(\vv)=0$.
\item[(2)] It satisfies Maxwell's equation
$$\nabla_\y\times\BS(\vv)(\y) = \vv(\y) + \nabla_\y\int_{\reals^3} \vv(\x)\cdot
\nabla_\x\ph_0(\x,\y)\,d\x,$$
where $\ph_0$ is the fundamental solution of the Laplacian in $\reals^3$.
\item[(3)] $\BS(\vv)(\y)\to \zero$ as $\y\to\infty$.
\end{enumerate}

To see that (2) is one of Maxwell's equations, first integrate by parts to rewrite
it as
$$\nabla_\y\times\BS(\vv)(\y) = \vv(\y) - \nabla_\y
\int_{\reals^3}(\nabla_\x\cdot\vv(\x))\,\ph_0(\x,\y)\,d\x.$$
If we think of the vector field $\vv(\x)$ as a steady current, then its negative divergence,
$-\nabla_\x\cdot\vv(\x)$, is the time rate of accumulation of charge at $\x$, and hence the
integral
$$-\nabla_\y\int_{\reals^3}(\nabla_\x\cdot\vv(\x))\,\ph_0(\x,\y)\,d\x$$
is the time rate of increase of the electric field $E$ at $\y$. Thus equation (2) is simply
Maxwell's equation
$$\nabla\times B = \vv +\frac{\partial E}{\partial t}.$$
In $\reals^3$, $S^3$ and $H^3$, a linear operator satisfying conditions (1), (2) and (3)
above will be referred to as a \it Biot-Savart operator\rm.

\vfill
\eject

\noindent\bf Remarks\rm.
\begin{itemize}
\item To see that equation (2) above is Maxwell's equation, we integrated
by parts, in spite of the fact that the kernel function $\ph_0(\alpha)$
has a singularity at $\alpha=0$. We leave it to the reader to check that
the validity of this depends on the fact that
the singularity of $\ph_0$ is of order $1/\alpha$. We will use this throughout
the paper, without further mention.
\item Recall Amp\`ere's Law: Given a divergence-free current flow, the
circulation of the resulting magnetic field around a loop is equal to the
flux of the current through any surface bounded by that loop. This is an
immediate consequence of Maxwell's equation (2) above, since if the current
flow $\vv$ is divergence-free, this equation says that $\nabla\times
\BS(\vv)=\vv$. Then Amp\`ere's Law is just the curl theorem of vector
calculus.

In particular, if the current flows along a wire loop, the circulation of
the resulting magnetic field around a second loop disjoint from it is
equal to the flux of the current through a cross-section of the wire loop,
multiplied by the linking number of the two loops. Thus linking numbers
are built into Amp\`ere's Law, and once we have an explicit integral
formula for the magnetic field due to a given current flow, we easily
get an explicit integral formula for the linking number.
\item In $\reals^3$, conditions (1), (2) and (3) are easily seen to
characterize the Biot-Savart operator, as follows. Since conditions
(1) and (2) specify the divergence and the curl of $\BS(\vv)$, the
difference $\BS_1(\vv)-\BS_2(\vv)$ between two candidates for the
Biot-Savart operator would be divergence-free and curl-free. Since
$\reals^3$ is simply connected, this difference would be the gradient
of a harmonic function. Hence the components of this gradient must also
be harmonic functions. Since they go to zero at infinity, they have
to be identically zero. Thus $\BS_1(\vv)=\BS_2(\vv)$.
\item In $S^3$, conditions (1) and (2) alone suffice to characterize
the Biot-Savart operator, since there are no non-zero vector fields on
$S^3$ which are simultaneously divergence-free and curl-free (i.e.,
there are no non-constant harmonic functions).
\item In $H^3$, it is not yet clear to us how to characterize the
Biot-Savart operator. Even strengthening (3) to require that
$\BS(\vv)(\y)$ go to zero exponentially fast at infinity is not
quite enough. And in $H^3$, unlike $\reals^3$, the field $\BS(\vv)$
is not in general of class $L^2$.
\end{itemize}
\vfill
\eject

The following theorem is from our first paper.

\medskip

\noindent\bf Theorem 2.1\rm. \sc Biot-Savart integrals in $S^3$ and
$H^3$\it. Biot-Savart operators exist in $S^3$ and $H^3$, and are
given by the following formulas, in which $\vv$ is a smooth, compactly
supported vector field:

\rm(1) \it On $S^3$, in left-translation format:
$$\eqalign{
\BS(\vv)(\y)=\int_{S^3}L_{\y\x^{-1}}\vv(\x)\times\nabla_\y
\ph_0&(\x,\y)\,d\x -\frac{1}{4\pi^2}\int_{S^3} L_{\y\x^{-1}}\vv(\x)
\,d\x\cr
&+2\,\nabla_\y\int_{S^3}L_{\y\x^{-1}}\vv(\x)\cdot\nabla_\y
\ph_1(\x,\y)\,d\x,}$$
\indent where $\ph_0(\alpha)=-(\pi-\alpha)\cot\alpha/(4\pi^2)$ and
$\ph_1(\alpha)=-\alpha(2\pi-\alpha)/(16\pi^2)$.

\rm(2) \it On $S^3$ in parallel transport format:
$$\BS(\vv)(\y)=\int_{S^3}P_{\y\x}\vv(\x)\times \nabla_\y\ph_0(\x,\y)\,d\x,$$
\indent where $\ph_0(\alpha)=-(\pi-\alpha)\csc\alpha/(4\pi^2)$.

\rm(3) \it On $H^3$ in parallel transport format:
$$\BS(\vv)(\y)=\int_{H^3}P_{\y\x}\vv(\x)\times \nabla_\y\ph_0(\x,\y)\,d\x,$$
\indent where $\ph_0(\alpha)=-\csch\alpha/(4\pi)$.\rm

\medskip

In formula (1), the function $\ph_1(\alpha)=-\alpha(2\pi-\alpha)/
(16\pi^2)$ satisfies the equation
$$\Delta\ph_1=\ph_0-[\ph_0],$$
where $[\ph_0]$ denotes the average value of $\ph_0$ over $S^3$.
The other kernel functions already appeared in the linking integrals
in Theorem 1.1.

In formula (3), the magnetic field $\BS(\vv)(\y)$ goes to zero at
infinity like $e^{-\alpha}$, where $\alpha$ is the distance from
$\y$ to a fixed point in $H^3$.

\vfill
\eject

\noindent\bf 3.\ Link, twist and writhe in $S^3$ and $H^3$\rm.

In a series of three papers (1959--1961), Georges C{\accent21 a}lug{\accent21 a}reanu defined
a real-valued invariant of a smooth \it simple \rm closed curve in $\reals^3$ by allowing  
the two curves in Gauss's linking integral to come together. 
In the limit, the points $\x(s)$ and $\y(t)$ now run along the same curve, and therefore can coincide,
making Gauss's integral seem improper because of the $|\x-\y|^3$ in the denominator.
But C{\accent21 a}lug{\accent21 a}reanu noted that in this case the numerator
goes to zero even faster than the denominator, so that the whole integrand goes to
zero as $\x$ and $\y$ come together, and the integral converges.
In (1971),
F.\ Brock Fuller called this invariant, which measures the extent to which the
curve wraps and coils around itself, the ``writhing number'':  
$$\Wr(K)  =  \int_{K\times K}
\frac{d\x}{ds}\times \frac{d\y}{dt}\cdot \frac{\x-\y}{4\pi|\x-\y|^3}\,ds\,dt.$$
In those papers, C{\accent21 a}lug{\accent21 a}reanu also discovered the formula  
\sc link = twist + writhe\rm, in which  \sc link \rm is the linking
number of the two edges of a closed ribbon, \sc twist \rm measures
the extent to which the ribbon twists around one of its edges, and
\sc writhe \rm is the writhing number of that edge.

C{\accent21 a}lug{\accent21 a}reanu proved this formula under the assumption that the
simple closed curve $K$ has nowhere-vanishing curvature, but the basic ideas for proving the formula without
this assumption are already present in his papers.  
This can be seen in sections 6 and 7 of this paper, where the proofs we give in $S^3$ and $H^3$ follow
C{\accent21 a}lug{\accent21 a}reanu's original proof in $\reals^3$, but require no curvature restriction.
Nevertheless, C{\accent21 a}lug{\accent21 a}reanu's
formula without the curvature restriction was proved by James White (1969)
in his thesis, using a totally different approach based on ideas of
William Pohl (1968a, b).

Moving on to $S^3$ and $H^3$, we follow C{\accent21 a}lug{\accent21 a}reanu's lead and replace the two 
closed curves $K_1$ and $K_2$ in the linking integrals of Theorem 1.1 by one simple closed curve.
Again all the integrals converge, and we use them to extend the notion of writhing number to these
spaces.

\vfill
\eject

\noindent\bf Definition of the writhing integrals in $S^3$ and
$H^3$\it.

\rm(1) \it On $S^3$ in left-translation format:
$$\eqalign{
\Wr_L(K) := \int_{K\times K} L_{\y\x^{-1}}\frac{d\x}{ds}&\times
\frac{d\y}{dt}\cdot \nabla_\y\ph(\x,\y)\,ds\,dt\cr
&\qquad -\frac{1}{4\pi^2}\int_{K\times K} L_{\y\x^{-1}}\frac{d\x}{ds}
\cdot \frac{d\y}{dt}\,ds\,dt,}
$$
\indent where $\ph(\alpha)= (\pi-\alpha)\cot\alpha/(4\pi^2)$.

\rm(2) \it On $S^3$ in parallel transport format:
$$\Wr_P(K) := \int_{K\times K} P_{\y\x}\frac{d\x}{ds}\times
\frac{d\y}{dt}\cdot \nabla_\y\ph(\x,\y)\,ds\,dt,$$
\indent where $\ph(\alpha)=(\pi-\alpha)\csc\alpha/(4\pi^2)$.

\rm(3) \it On $H^3$ in parallel transport format:
$$\Wr_P(K) := \int_{K\times K} P_{\y\x}\frac{d\x}{ds}\times
\frac{d\y}{dt}\cdot \nabla_\y\ph(\x,\y)\,ds\,dt,$$
\indent where $\ph(\alpha)=\csch\alpha/(4\pi)$.\rm

\medskip

The two versions of the writhing number on $S^3$ are \it not \rm the
same, and one can show that
$$\Wr_L(K)=\Wr_P(K)+\frac{\mbox{\rm length of $K$}}{2\pi}.$$
The parallel transport
version of writhe is more intuitively satisfying, since in this version
the writhing number of a great circle is zero.

We turn next to the definition of ``twist''.

Let $K$ be a smooth simple closed curve in $S^3$ or $H^3$, parametrized
by arclength $s$. Let $\x(s)$ be a moving point
along $K$, and $\T(s)=\x'(s)$ the unit tangent vector field.

Let $\vv(s)$ be a unit normal vector field along $K$. Our
intention is to define the (total) twist of $\vv$ along $K$ by a
formula similar to
$$\Tw(\vv)=\frac{1}{2\pi}\int_K \T(s)\times\vv(s)\cdot\vv'(s)\,ds,$$
the formula for twist in Euclidean 3-space.

%%% FIGURE %%%
\begin{figure}[h!]
\center{\includegraphics[height=130pt]{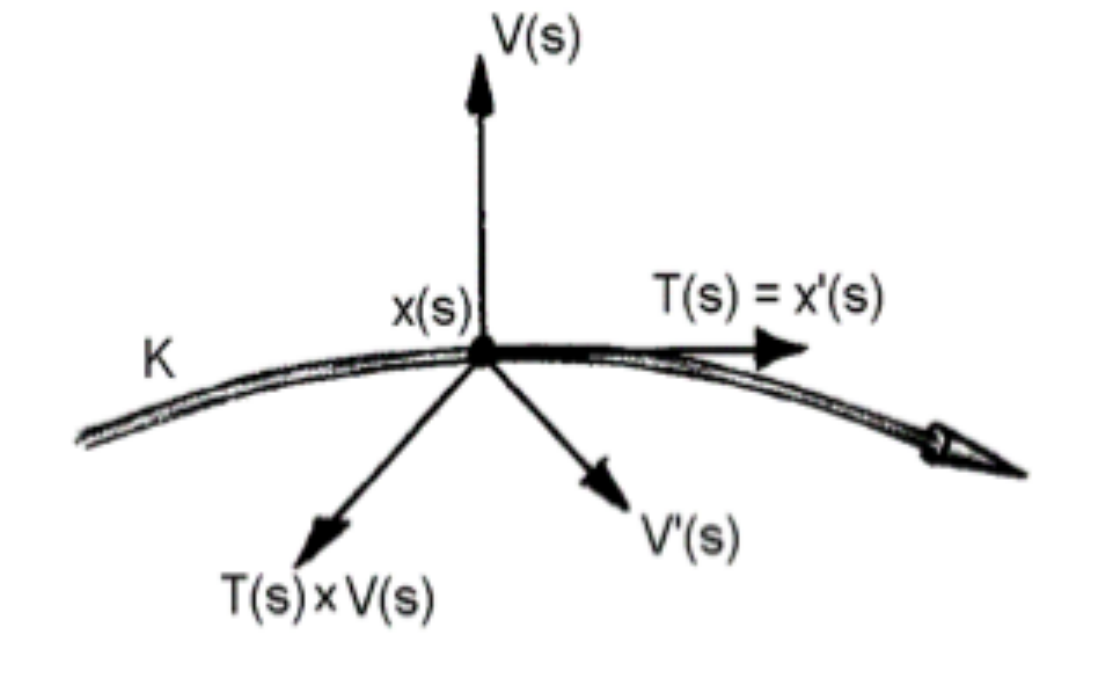}
\caption{Vectors in the definition of twist}}
\end{figure}
%%%%%%%%%%%%%%%%%%%%

But on $S^3$ there are two flavors of twist, according as $\vv'(s)$
is calculated as a ``left-invariant'' derivative or as a covariant
derivative. If we fall back into Euclidean mode and write
$$\vv'(s)=\lim_{h\to 0}\frac{\vv(s+h)-\vv(s)}{h},$$
then the vectors $\vv(s+h)$ and $\vv(s)$ lie in different tangent spaces,
and we must move them together in order to subtract. If we use
left-translation in the group $S^3$ to move $\vv(s+h)$ back to
the tangent space at $\x(s)$ which contains $\vv(s)$, then the
resulting limit is the \it left-invariant derivative \rm $\vv'_L(s)$.
If we use parallel transport to move $\vv(s+h)$ back, then the resulting
limit is the \it covariant derivative \rm $\vv'_P(s)$.

The two flavors of twist on $S^3$ are then given by
$$\Tw_L(\vv)=\frac{1}{2\pi}\int_K\T(s)\times\vv(s)\cdot\vv'_L(s)\,ds
\eqno(1)$$
and
$$\Tw_P(\vv)=\frac{1}{2\pi}\int_K\T(s)\times\vv(s)\cdot\vv'_P(s)\,ds.
\eqno(2)$$
One can show that
$$\Tw_L(\vv)=\Tw_P(\vv)-\frac{\mbox{\rm length of $K$}}{2\pi}.$$

\medskip

\noindent\bf Example\rm. Consider the great circle $K =\{(\cos s,\,
\sin s,\,0,\,0)\,:\,0\le s\le 2\pi\}$ on $S^3$, and along it
the unit normal vector field $\vv(s)=(0,0,\cos s,\sin s)$. Then we
have $\Tw_L(\vv)=0$ and $\Tw_P(\vv)=1$.

\medskip

In hyperbolic 3-space $H^3$, we have only the parallel transport version
of twist,
$$\Tw_P(\vv)=\frac{1}{2\pi}\int_K\T(s)\times\vv(s)\cdot\vv'_P(s)\,ds.
\eqno(3)$$

Now consider in $S^3$ or $H^3$ a narrow ribbon of width $\eps$ obtained
by starting with a simple closed curve $K=\{\x(s)\}$ and then
exponentiating a unit normal vector field $\vv(s)$ along $K$. One
edge of this ribbon is the original curve $K$, and the other edge is
the curve $K_\eps=\{\y_\eps(s)\}$, given explicitly (see section 5) by
$$\y_\eps(s)=\cos\eps\ \x(s)+\sin\eps\ \vv(s)\ \ \mbox{\rm in\ }
\ S^3;$$
$$\y_\eps(s)=\cosh\eps\ \x(s)+\sinh\eps\ \vv(s)\ \ \mbox{\rm in\ }
\ H^3.$$
Since $K$ is simple, the ribbon will be \it embedded \rm in $S^3$ or
$H^3$ provided $\eps$ is small enough.

The vector field $\vv(s)$ then points ``across'' the ribbon.

%%% FIGURE 3 %%%
\begin{figure}[h!]
\center{\includegraphics[height=160pt]{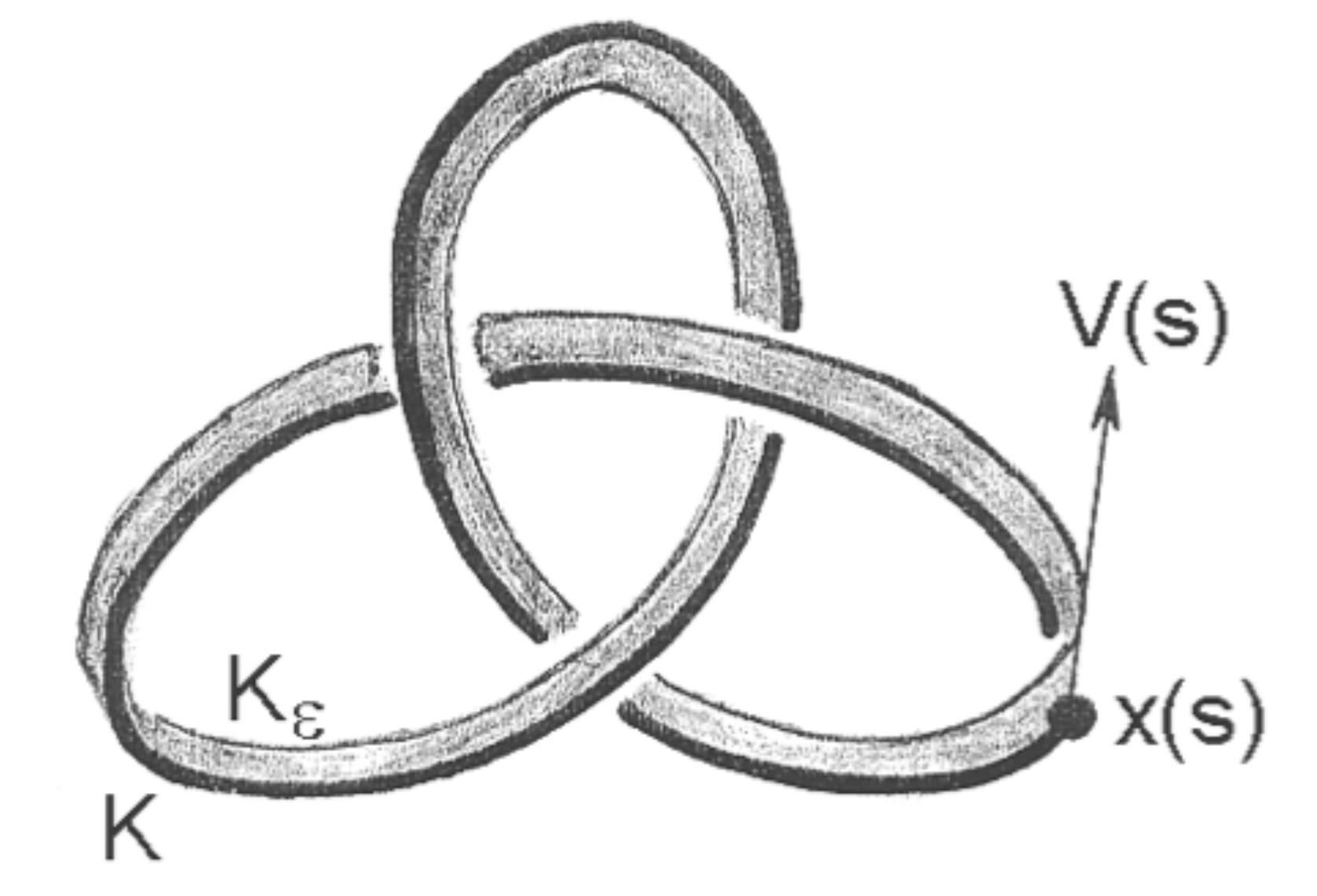}
\caption{A ribbon, its generating curve, and its vector field}}
\end{figure}
%%%%%%%%%%%%%%%%%%%%

\noindent\bf Theorem 3.1\rm. \sc link = twist + writhe \rm in $S^3$ and
$H^3$.

\rm(1) \it On $S^3$ in left-translation format:
$$\Lk(K,K_\eps)= \Tw_L(\vv)+\Wr_L(K).$$

\rm(2) \it On $S^3$ in parallel transport format:
$$\Lk(K,K_\eps)= \Tw_P(\vv)+\Wr_P(K).$$

\rm(3) \it On $H^3$ in parallel transport format:
$$\Lk(K,K_\eps)= \Tw_P(\vv)+\Wr_P(K).$$\rm

\noindent We give an overview of the proof in the next section.

\vfill
\eject
\noindent\bf 4.\ Proof scheme for \sc link = twist + writhe\rm.

In spirit, our proof of Theorem 3.1 for ribbons in $S^3$ and $H^3$
follows C{\accent21 a}lug{\accent21 a}reanu's original proof in $\reals^3$: we begin with
the linking integrals given in Theorem 1.1 for the edges $K$ and $K_\eps$
of our ribbon, let $\eps$ shrink to zero, and observe the behavior
of the linking integrand.

The value of the linking integral is independent of $\eps$
for $\eps>0$ since the ribbon is embedded and since the linking number 
is invariant under homotopies which keep the two curves disjoint. 
But the linking integrand blows up as one approaches the
diagonal of $K\times K$, and this is
handled as follows.

Outside an appropriately chosen neighborhood of the diagonal, the
linking integrand converges to the writhing integrand as $\eps\to 0$,
and its integral converges to the writhing number of the curve $K$.
Inside this neighborhood of the diagonal, the linking integrand
blows up, but its integral converges to the total twist of the
normal vector field $\vv$ along $K$.

The crucial thing, recognized by C{\accent21 a}lug{\accent21 a}reanu, is that the width of
the neighborhood of the diagonal in $K\times K$ must go to zero
\it much more slowly \rm than the width $\eps$ of the ribbon. In fact,
we will choose the neighborhood of the diagonal to have width
$\eps^p$, where $0<p<1/3$.

To give a sense of this in action, we will outline here the proof of
Theorem 3.1, part (2), dealing with \sc link = twist + writhe \rm in
parallel transport format on $S^3$. The proofs for $H^3$ and for
left-translation format on $S^3$ are essentially the same. In
particular, in left-translation format, the integrand of the second
integral in the expression for the linking number converges uniformly
to the corresponding integrand for the writhing number.

Consider, in parallel transport format on $S^3$, the linking integrand
of $K$ with $K_\eps$,
$$F_\eps(s,t)=\frac{d\x}{ds}\cdot
P_{\x(s)\y_\eps(t)}\left(\frac{d\y_\eps}{dt}\times
\nabla_{\y_\eps(t)}\ph\bigl(\alpha(\x(s),\y_\eps(t))\bigr)\right)$$
and the writhing integrand of $K$,
$$F_0(s,t)=\frac{d\x}{ds}\cdot
P_{\x(s)\x(t)}\left(\frac{d\x}{dt}\times
\nabla_{\x(t)}\ph\bigl(\alpha(\x(s),\x(t))\bigr)\right),$$
where $\ph(\alpha)=(\pi-\alpha)\csc\alpha/(4\pi^2)$.

\vfill
\eject
Then the linking number of $K$ and $K_\eps$ is given by
$$\Lk(K,K_\eps)=\int\!\!\int_{0\le s,t\le L}F_\eps(s,t)\,ds\,dt,$$
and the writhing number of $K$ is given by
$$\Wr_P(K)=\int\!\!\int_{0\le s,t\le L}F_0(s,t)\,ds\,dt.$$
Since $\alpha(\x,\y)$ is the distance between $\x$ and $\y$, and
since $\ph$ has a singularity just at $\alpha=0$, the only difficulty
in considering the convergence of the linking integral as $\eps\to 0$
happens near the diagonal, where $s=t$.

To handle this, we first show that because the singularity of
$\ph''(\alpha)$ at $\alpha=0$ is like $1/\alpha^3$, we have that
for sufficiently small $\eps>0$,
$$|F_\eps(s,t)-F_0(s,t)|\le C\eps^{1-3p},$$
provided that $|s-t|\ge\eps^p$.

If $0<p<1/3$, then $1-3p>0$, and hence $C\eps^{1-3p}\to 0$ as
$\eps\to 0$. Therefore $F_\eps(s,t)$ converges uniformly to
$F_0(s,t)$ in the region $|s-t|\ge\eps^p$, and this region expands
to the region $|s-t|>0$ as $\eps\to 0$. Since the writhing
integrand $F_0(s,t)$ remains bounded even along the diagonal, this
shows that
$$\int\!\!\int_{|s-t|>\eps^p}F_\eps(s,t)\,ds\,dt \quad \xrightarrow{\hphantom{=}\eps\to 0
\hphantom{=}}
\quad\int\!\!\int_{0\le s,t\le L}F_0(s,t)\,ds\,dt=\Wr_P(K),$$
that is,  a \it portion \rm of the linking integral
converges to the \it entire \rm writhing integral as $\eps\to 0$. This is the content
of Proposition 6.3 below.

The more delicate part of the argument is the integral near the
diagonal. A careful analysis reveals that
for $0<p<1$,
$$\lim_{\eps\to 0}\int_{t-\eps^p}^{t+\eps^p}F_\eps(s,t)\,ds =
\frac{1}{2\pi}\x'(t)\times\vv(t)\cdot \vv'_P(t).$$
Hence
$$\int\!\!\int_{|s-t|<\eps^p}F_\eps(s,t)\,ds\,dt\quad\xrightarrow{\hphantom{=}\eps\to 0
\hphantom{=}}
\quad\frac{1}{2\pi}\int_0^L\x'(t)\times\vv(t)\cdot
\vv'_P(t)\,dt=\Tw_P(\vv).$$
That is, the remaining portion of the linking integral converges to
the entire twisting integral. This is the content of Proposition 6.4. In this way, we see that
$$\Lk(K, K_\eps)=\Tw_P(\vv)+\Wr_P(K).$$

\bigskip

\noindent \bf 5.\ Some geometric formulas on $S^3$\rm

Before we can proceed with the details of the proof of \sc link
= twist + writhe\rm, we need to collect some basic geometric
formulas on $S^3$, which are treated in more detail in our (2008) paper.

We consider $S^3\subset\reals^4$ in the usual way, as the set
$$\{\x\in\reals^4\,|\,\ip{\x}{\x}=1\},$$
where $\ip{\x}{\y}$ is the standard inner product on $\reals^4$.
Since the linking, twisting and writhing integrands involve
cross-products of vectors, we remind the reader that
if \\
$\x\in S^3$, and $\vv,\w\in T_\x S^3$, we define the cross product by
$$\vv\times\w = \det\left[ \begin{array}{cccc}
x_0 & x_1 & x_2 & x_3 \\
v_0 & v_1 & v_2 & v_3 \\
w_0 & w_1 & w_2 & w_3 \\
\xh_0 & \xh_1 & \xh_2 & \xh_3 \end{array}
\right].$$
In this formula, we view $\x$, $\vv$, $\w$ and the result as vectors in $\reals^4$ and
$\{\xh_0,\xh_1,\xh_2,\xh_3\}$ is the canonical orthonormal basis of
$\reals^4$. From this, it is easy to see that if $\u$ is also tangent to
$S^3$ at $\x$, then the triple product $\u\cdot\vv\times\w$ is equal to
the value of the 4-by-4 determinant whose rows are $\x$, $\u$, $\vv$ and
$\w$. We will use the notation $\Vert \x,\u,\vv,\w\Vert$ for this
determinant.

Next, suppose $\vv$ is a \it unit \rm vector in $T_\x S^3$. Then the unique
unit-speed geodesic in $S^3$ through $\x$ with initial
tangent vector $\vv$ is given by
$$G(t)=\cos t\ \x + \sin t\ \vv.$$
Because $\ip{\x}{\vv}=0$, we have that
$\ip{\x}{G(t)}=\cos t$, and we can conclude in general that the
geodesic distance $\alpha(\x,\y)$ between two points $\x$ and $\y$
on $S^3$ is
$\alpha(\x,\y)=\arccos\ip{\x}{\y}$.

Moreover, if $\x$ and $\y$ are any distinct, non-antipodal points on
$S^3$,
then the vector $\vv=(\y-\cos\alpha\ \x)/\sin\alpha$ is a unit vector in
$T_\x S^3$, and the geodesic
it generates connects $\x$ to $\y$. From this we deduce that
$$\nabla_\x\alpha(\x,\y)=\frac{\cos\alpha\ \x-\y}{\sin\alpha}$$
and
$$\nabla_\y\alpha(\x,\y)=\frac{\cos\alpha\ \y-\x}{\sin\alpha}.$$

We will also need the formula for parallel transport of a vector
$\vv\in T_\y S^3$ to $T_\x S^3$:
$$P_{\x\y}(\vv) = \vv - \frac{\ip{\x}{\vv}}{1+\ip{\x}{\y}}(\x+\y).$$
Specifically, we need the observation that $P_{\x\y}$ affects $\vv$
by adding a linear combination of $\x$ and $\y$.

Using these formulas, we can make precise the definitions of $F_\eps(s,t)$ and
$F_0(s,t)$\\and then derive equivalent expressions for them that will be useful in
our proof of\\ 
\sc link = twist + writhe\rm. 

\medskip

\noindent\bf Proposition 5.1\it.\ \ Let $\x(s)$ be a simple closed
curve in $S^3$, let $\vv(s)$ be a
unit vector in $T_{\x(s)}S^3$ which is perpendicular to $\x'(s)$, and set
$\ye(t)=\x(t)\cos\eps+\vv(t)\sin\eps$ for each $t$. Then
$$\eqalign{
F_\eps(s,t)& =
\frac{d\x}{ds}\cdot P_{\x(s)\ye(t)}\left(\frac{d\ye}{dt}\times\nabla_{\ye(t)}
\ph\bigl(\alpha(\x(s),\ye(t))\bigl)\right) \cr
&=-\frac{\ph'(\alpha_\eps)}{\sin\alpha_\eps} \left\Vert\ye(t),\frac{d\ye}{dt},
\x(s),\frac{d\x}{ds}\right\Vert
}$$
using the determinant notation given above and the shorthand 
$\aes=\alpha(\x(s),\ye(t))$ for the distance between $\x(s)$ and $\ye(t)$.
Similarly, we have (for $s\ne t$)
$$\eqalign{
F_0(s,t)&=\frac{d\x}{ds}\cdot
P_{\x(s)\x(t)}\left(\frac{d\x}{dt}\times
\nabla_{\x(t)}\ph\bigl(\alpha(\x(s),\x(t))\bigr)\right)\cr
&= -\frac{\ph'(\alpha_0)}{\sin\alpha_0} \left\Vert\x(t),\frac{d\x}{dt},\x(s),
\frac{d\x}{ds}\right\Vert
}$$
where $\alpha_0$ is the distance between $\x(s)$ and $\x(t)$. \rm

%%% FIGURE 3 %%%
\begin{figure}[h!]
\center{\includegraphics[height=180pt]{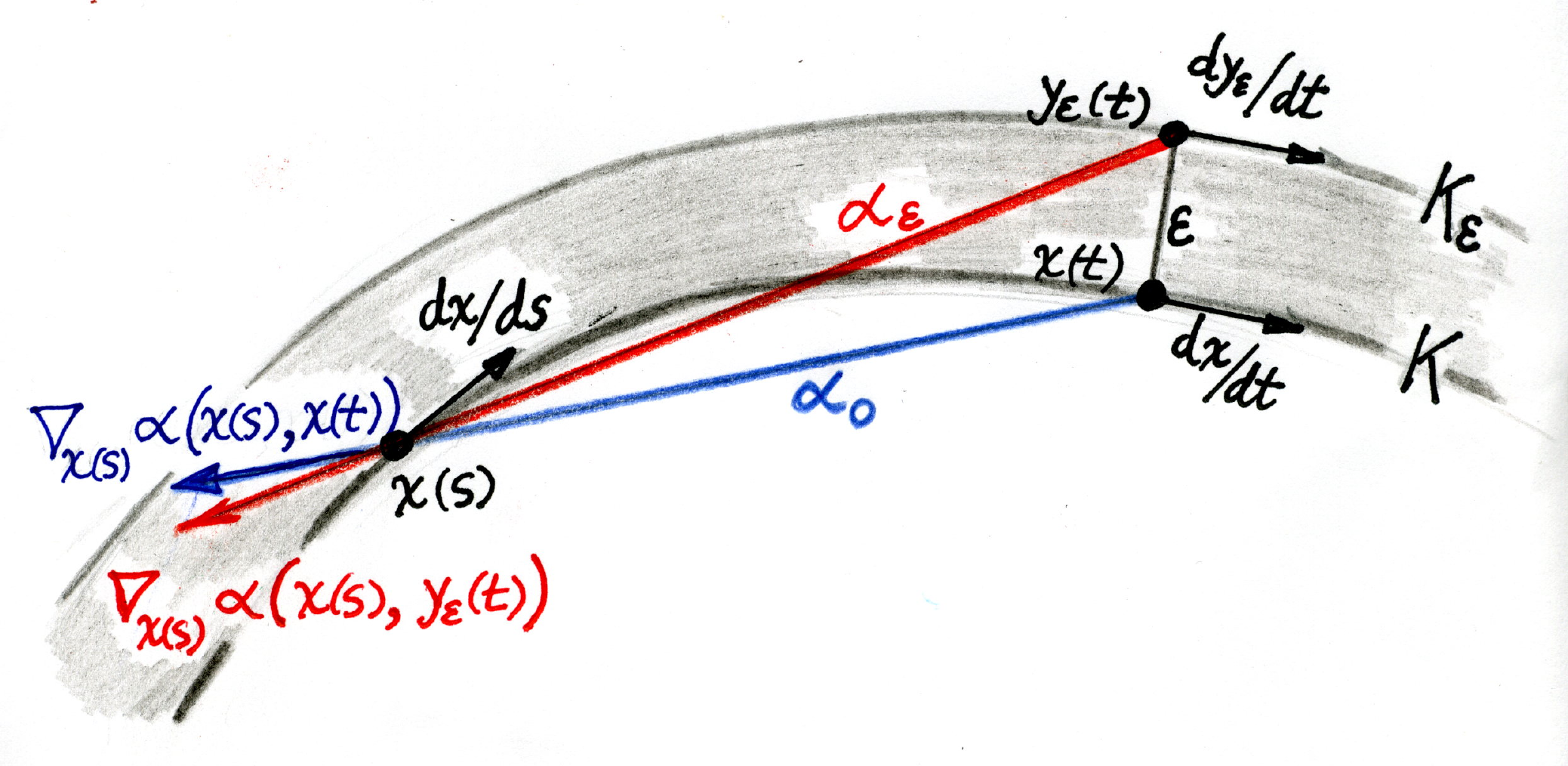}
\caption{A ribbon, its generating curve, and relevant vectors}}
\end{figure}
%%%%%%%%%%%%%%%%%%%%

\noindent\it Proof\rm.\ \ Using the formulas given above for $\nabla_{\ye}
\alpha$ and the cross product, we write
$$\eqalign{
\frac{d\ye}{dt}\times\nabla_{\ye} \ph(\alpha(\x,\ye(t)))
&=\ph'(\alpha_\eps)\frac{d\ye}{dt}\times\frac{\cos\alpha_\eps\,\ye-\x}{\sin
\alpha_\eps}\cr
&= -\frac{\ph'(\alpha_\eps)}{\sin\alpha_\eps} \left[\ye(t),\frac{d\ye}{dt},\x(s)\right].}$$
Because the triple product is perpendicular to $\x$ and $\ye$,
this vector in $\reals^4$ is not changed by $P_{\x\ye}$.
Therefore we can express
$$F_\eps(s,t)=\frac{d\x}{ds}\cdot\left(
-\frac{\ph'(\alpha_\eps)}{\sin\alpha_\eps} \left[\ye(t),\frac{d\ye}{dt},\x(s)\right]\right)=
-\frac{\ph'(\alpha_\eps)}{\sin\alpha_\eps} \left\Vert\ye(t),\frac{d\ye}{dt},
\x(s),\frac{d\x}{ds}\right\Vert.$$
The proof for $F_0$ is identical.

\bigskip

\noindent \bf 6.\ Proof of \sc link = twist + writhe \bf on $S^3$
\bf in parallel transport format\rm

In this section we prove the \sc link = twist + writhe \rm
formula in parallel transport format for ribbons
in the 3-sphere. As outlined above, the
idea is to write the linking integral for the two edges of a ribbon
of width $\eps$, and then take its limit as $\eps\to 0$. Of course
the value of the linking integral stays constant, but the limit of
the integral is not equal to the integral of the limit of the linking
integrand. The latter limit is the writhe of the fixed
edge of the ribbon, and the difference is the twist.

To avoid unnecessary complications, we assume all
our curves and deformations of curves are smooth, so we are free to
differentiate, commute derivatives, etc.

As indicated above, we begin with a smooth, simple closed curve $K$ 
parametrized by
arclength and given by $\x(s)$ for $0\le s\le L = $ length of $K$. We
define our ribbon by letting $\vv(s)$ be a \it unit \rm vector, tangent to
$S^3$ at $\x(s)$ and perpendicular to
$\T(s)=\x'(s)$ for every $s$. The other edge of our ribbon of width
$\eps$ will be at distance
$\eps$ along the geodesic emanating from $\x$ in the
direction of $\vv$, so it is $\ye(s) =\cos\eps\ \x(s) +\sin\eps\ \vv(s)$.
In general, $s$ is not an arclength parameter for the curve $\ye(s)$.

The linking number of the two edges of the ribbon is:
$$\Lk(\x,\ye)=\int_0^L\!\!\!\int_0^L F_\eps(s,t)\,ds\,dt,$$
where $F_\eps(s,t)$ is given by either of the expressions in 
Proposition 5.1.

The linking number is
independent of $\eps$, and so our
strategy will be to take the limit of the linking integral of the
edges of the ribbon of
width $\eps$ as $\eps\to 0$. We will examine the
difference between the limit of the integral (the linking
number) and the integral of the limit (the writhing number), and
show that it is equal to the twist of the ribbon as defined earlier.
Since the twist is defined by a single integral
in contrast to the double integrals that define link and writhe, we'll
use the following notation for the ``halfway'' integrations of the
latter two quantities:
$$H_{Lk}(t;\eps)=
\int_0^L F_\eps(s,t)\,ds$$
and
$$H_{Wr}(t)=\int_0^L F_0(s,t)\,ds.$$
Our objective will be to show that
$$\lim_{\eps\to 0}H_{Lk}(t;\eps)=H_{Wr}(t)
+ \mbox{\rm``something''},$$
where the integral of ``something'' with respect to $t$
will be the twist of the ribbon.

As we indicated above, the convergence of the
linking integrand to the writhing integrand fails to
be uniform only near the diagonal of $[0,L]\times[0,L]$, so we'll write
$$H_{Lk}(t)=\int_{|s-t|>\eps^p}\cdots\,ds +
\int_{|s-t|<\eps^p}\cdots\,ds$$
where $p$ is a number between 0 and 1 to be determined later.
We will show that the first term converges to $H_{Wr}(t)$ and the second
term will give us our ``something''.

Before we can prove that the convergence is uniform away from the
diagonal, we need the following preliminary lemma.

\medskip

\noindent\bf Lemma 6.1\rm:\ \ \it There is a constant $C>0$
such that $\alpha(\x(s),\x(t))>C|s-t|$,
where we consider $|s-t|$ to be the
``distance'' on the circle with circumference $L$.

\noindent Proof\rm.\ \ This is true locally
(i.e., for $s$ near $t$) because $\x$ is parametrized by arclength
(as we will justify below), and globally by compactness.

To get the local estimate, we use Taylor's formula to write
$$\x(s)=\x(t)+(s-t)\x'(t)+\frac{(s-t)^2}{2}\x''(t)+(s-t)^3\A(s,t)$$
where $\A(s,t)$ is a bounded, smooth vector-valued function of $s$ and $t$.
Since $\x(s)$ lies on the sphere $S^3\subset \reals^4$, we have $\ip{\x(s)}{\x'(s)}=0$,
and since $\x$ is parametrized by arclength, we have $\ip{\x'(s)}{\x'(s)}=1$.
It follows that $\ip{\x(s)}{\x''(s)}=-1$, and hence  
$$\cos\alpha(\x(s),\x(t))=\ip{\x(s)}{\x(t)}=
1-\frac{(s-t)^2}{2}+(s-t)^3p_1(s,t),$$
where $p_1(s,t)$ is a bounded smooth scalar-valued function of
$s$ and $t$. In what
follows, $p_i$ will always stand for such a function without comment.
Then clearly 
$$\sin^2(\alpha)=1-\cos^2(\alpha)=(s-t)^2(1+(s-t)p_2)$$ 
and using
Taylor's theorem for $(1+z)^{1/2}$ and for $\arcsin(z)$ we conclude that
$$\alpha(\x(s),\x(t))=|(s-t)+(s-t)^2p_3|.$$
This is surely larger than $\frac{1}{2}|s-t|$ for $|s-t|$ sufficiently
small, say for $|s-t|<\delta$.

\noindent\bf Corollary 6.2\rm: \ \ \it
$\alpha(\x(s),\ye(t))>C'|s-t|$, with
$C'$ independent of $\eps$, provided
$\eps$ is small enough so that the ribbon never touches itself. When
$|s-t|>\eps^p$, this implies $\alpha(\x(s),\y_\eps(t))>C'\eps^p$.\rm

Again, this is a combination of a local estimate and a global compactness
argument.

\medskip

Now we can begin to analyze the convergence of the linking integral. We start with
the part away from the diagonal, which we expect to converge to the writhing
integral.

\medskip

\noindent\bf Proposition 6.3\rm:\ \ \it If $0<p<1/3$, then
$$
\lim_{\eps\to 0}
\int\!\!\!\int_{|s-t|>\eps^p} F_\eps(s,t)\,ds\,dt =
\int_0^L\int_0^L F_0(s,t)\,ds\,dt = \Wr(K),$$
in other words, the limit of the ``away from the diagonal'' part of
$\Lk(\x,\ye)$ is the integral of
$H_{Wr}(\x)$, which is the writhing number $\Wr(K)$.\rm

\medskip

\noindent\it Proof\rm.\ \ We need to
analyze the
difference $$F_\eps(s,t)-F_0(s,t) =
-\frac{\ph'(\alpha_\eps)}{\sin\alpha_\eps} \left\Vert\ye(t),\frac{d\ye}{dt},
\x(s),\frac{d\x}{ds}\right\Vert+\frac{\ph'(\alpha_0)}{\sin\alpha_0} \left\Vert\x(t),\frac{d\x}{dt},\x(s),
\frac{d\x}{ds}\right\Vert,$$
using the notation of Proposition 5.1.
 
Using properties of the determinant, we can rewrite the difference $F_\eps-F_0$
as a sum as follows:
$$F_\eps-F_0=\frac{-\ph'(\alpha_\eps)}{\sin\alpha_\eps}\left\Vert\ye(t),
\frac{d\ye}{dt}-\frac{d\x}{dt},
\x(s),\frac{d\x}{ds}\right\Vert+
\left\Vert\frac{\ph'(\alpha_0)}{\sin\alpha_0}\x(t)-
\frac{\ph'(\alpha_\eps)}{\sin\alpha_\eps}\ye(t),\frac{d\x}{dt},\x(s),\frac{d\x}{ds}\right\Vert.$$
We proceed to bound these two summands in terms of $\eps$.

For the first summand of $F_\eps-F_0$,
$$\frac{-\ph'(\alpha_\eps)}{\sin\alpha_\eps}\left\Vert\ye(t),\frac{d\ye}{dt}-\frac{d\x}{dt},
\x(s),\frac{d\x}{ds}\right\Vert,$$
we begin with some easy preliminary observations: $\x(s)$ and $d\x/ds$ are
unit vectors, and since $\ye(t)=\cos\eps\ \x(t)+\sin\eps\ \vv(t)$, we have
$$\frac{d\ye}{dt}-\frac{d\x}{dt}= (\cos\eps-1)\ \frac{d\x}{dt}+
\sin\eps\ \frac{d\vv}{dt},$$
and so can bound the second vector in the determinant as
$$\left|\frac{d\ye}{dt}-\frac{d\x}{dt}\right|\le C\eps,$$
where $C$ depends on the maximum value of $|d\vv/dt|$.

To handle the first vector in the determinant, we'll group the $1/\sin\alpha_\eps$
with the $\ye(t)$, and then note that the determinant is unaffected if we
subtract $(\cos\alpha_\eps/\sin\alpha_\eps)\x(s)$ from the first vector. In other
words, the first summand of $F_\eps-F_0$ is equal to
$$-\ph'(\alpha_\eps)\left\Vert\frac{1}{\sin\alpha_\eps}\ye(t)-
\frac{\cos\alpha_\eps}{\sin\alpha_\eps}\x(s),\frac{d\ye}{dt}-\frac{d\x}{dt},\x(s),
\frac{d\x}{ds}\right\Vert.$$
And since $\cos\alpha_\eps=\ip{\x(s)}{\ye(t)}$, the first vector in this latter
determinant is a unit vector.
Therefore, the entire determinant is bounded by $|\ph'(\alpha_\eps)|C\eps$.

Finally, recall from section 1 that we can bound  $|\ph'(\alpha)|$
by a constant divided by $\alpha^2$,
and since
$\alpha_\eps>\eps^p$ by hypothesis, we conclude that the first summand
of $F_\eps-F_0$ is bounded by
$M\eps^{1-2p}$.

The second summand of $F_\eps-F_0$ is the determinant
$$\left\Vert \frac{\ph'(\alpha_0)}{\sin\alpha_0}\,\x(t)-
\frac{\ph'(\alpha_\eps)}{\sin\alpha_\eps}\,\ye(t),\frac{d\x}{dt},\x(s),\frac{d\x}{ds}\right\Vert,$$
and our job will be to handle its first vector, since the other three are all
unit vectors.

Since the value of the determinant 
is unaffected if we replace its first vector
with
$$\ph'(\alpha_0)\left(\frac{1}{\sin\alpha_0}\,\x(t)-\frac{\cos\alpha_0}{\sin\alpha_0}\,\x(s)\right)-\ph'(\alpha_\eps)\left(\frac{1}{\sin\alpha_\eps}
\,\ye(t)-\frac{\cos\alpha_\eps}{\sin\alpha_\eps}\,\x(s)\right),$$
we will obtain a bound on the determinant by bounding this vector.
Using the expressions for $\nabla\alpha$ derived in section 5,
we write this as
$$\nabla_{\x(s)}\ph\bigl(\alpha(\x(s),\ye(t))\bigr)-\nabla_{\x(s)}
\ph\bigl(\alpha(\x(s),\x(t))\bigr),$$
and then rewrite it as
$$\int_0^\eps \frac{d}{d\sigma}
\nabla_{\x(s)}\ph\bigl(\alpha(\x(s),\ys(t))\bigr)\,d\sigma,$$
where $\ys(t)=\cos\sigma\,\x(t)+\sin\sigma\,\vv(t)$.

We now calculate and estimate:
$$\eqalign{
\frac{d}{d\sigma}
\nabla_{\x(s)}\ph\bigl(\alpha(\x(s),\ys(t))\bigr)&=
\frac{d}{d\sigma}\left(\ph'\bigl(\alpha(\x(s),\ys(t))\bigr)\nabla_{\x(s)}
\alpha(\x(s),\ys(t))\right)\cr
&=I + II,}$$
where $$I=\ph''\bigl(\alpha(\x(s),\ys(t))\bigr)\frac{d\alpha}{d\sigma}\ \nabla_{\x(s)}
\alpha(\x(s),\ys(t))$$
and
$$II=\ph'\bigl(\alpha(\x(s),\ys(t))\bigr)\,\frac{d}{d\sigma}
\nabla_{\x(s)}\alpha(\x(s),\ys(t)).$$

To bound $|I|$, we know that $|\nabla\alpha|=1$,
$|d\alpha/d\sigma|\le 1$,
and $\ph''(\alpha)=1/(2\pi\alpha^3)+\cdots$. Since we also know that
$\alpha>K'\eps^p>K'\sigma^p$, we get 
$$|I|<\frac{Q_1}{\sigma^{3p}}.$$

To bound $|II|$, we have to know more about
$$\eqalign{
\frac{d}{d\sigma}\nabla_{\x(s)}\alpha(\x(s),\ys(t))&=
\frac{d}{d\sigma}\left(\frac{\cos\as\,\x-\ys}{\sin\as}\right)
\cr
&= -\frac{1}{\sin^2\as}
\frac{d\alpha}{d\sigma}\x-\frac{1}{\sin\as}\,\frac{d\ys}{d\sigma}+\frac{\cos\as}{\sin^2\as}
\frac{d\alpha}{d\sigma}\ys \cr
&= -\frac{1}{\sin\as}\frac{d\ys}{d\sigma}+\frac{1}{\sin\as}\frac{d\alpha}{d\sigma}
\left(\frac{\cos\as\,\ys-\x}{\sin\as}\right)\cr
&=-\frac{1}{\sin\as}\frac{d\ys}{d\sigma}+\frac{1}{\sin\as}
\frac{d\alpha}{d\sigma}\nabla_{\ys}\alpha.}$$
Once again, we'll use the facts that $|\nabla\alpha|=1$,
$|d\alpha/d\sigma|\le 1$, and $|d\ys/d\sigma|=1$ to conclude that
$$\left|\frac{d}{d\sigma}\nabla_\x \alpha(\x,\ys)\right|\le
\left|\frac{2}{\sin\alpha}\right|<\frac{C}{\alpha}.$$
Finally, we use that $\ph'(\alpha)=-1/(4\pi\alpha^2)+\cdots$, so that
$|\ph'(\alpha)|\le C'/\alpha^2$, and the fact that
$\alpha>K'\eps^p>K'\sigma^p$
to conclude that
$$|II|\le|\ph'|\left|\frac{d}{d\sigma}\nabla\alpha\right|\le \frac{Q_2}{\sigma^{3p}}.$$

Now we've estimated both terms into which we decomposed
$(d/d\sigma)
\nabla_{\x(s)}\ph\bigl(\alpha(\x(s),\ys(t))\bigr)$, so we can estimate its integral
as $\sigma$ goes from 0 to $\eps$ to obtain the result
$$\left\vert \nabla_{\x(s)}\ph\bigl(\alpha(\x(s),\ye(t))\bigr)-\nabla_{\x(s)}
\ph\bigl(\alpha(\x(s),\x(t))\bigr)\right\vert
\le\int_0^\eps \frac{Q_1+Q_2}{\sigma^{3p}}\,d\sigma
=Q\eps^{1-3p}.$$

So far, for $|s-t|>\eps^p$, we have
$$|F_\eps(s,t)-F_0(s,t)|\le M\eps^{1-2p}+Q\eps^{1-3p}.$$
Therefore
$$\left|\int_{|s-t|>\eps^p}F_\eps(s,t)\,ds-\int_{|s-t|>\eps^p}F_0(s,t)\,ds\right|
\le L(M\eps^{1-2p}+Q\eps^{1-3p}).$$
Since
$$\left|\int_{t-\eps^p}^{t+\eps^p}F_0(s,t)\,ds\right|\le R\eps^p,$$
we get
$$\left|\int_{|s-t|>\eps^p}F_\eps(s,t)\,ds-\int_0^L F_0(s,t)\,ds\right|
\le LM\eps^{1-2p}+LQ\eps^{1-3p}+R\eps^p.$$
So if $0<p<1/3$, we can conclude that
$$\lim_{\eps\to 0}
\int_{|s-t|>\eps^p} F_\eps(s,t)\,ds =
H_{Wr}(t)$$
uniformly in $t$, and so
$$
\lim_{\eps\to 0}
\int\!\!\!\int_{|s-t|>\eps^p} F_\eps(s,t)\,ds\,dt =
\int_0^L\!\!\int_0^L F_0(s,t)\,ds\,dt = \Wr(K).$$
This completes the proof of Proposition 6.3.

\medskip

Now we must analyze the part of the linking integral near the
diagonal.

\medskip

\noindent\bf Proposition 6.4\it:\ \ With $\x$, $\ye$ and $F_\eps$
defined as above, for $0<p<1/3$,
$$\lim_{\eps\to 0}\int_{t-\eps^p}^{t+\eps^p}F_\eps(s,t)
 \,ds =
\frac{1}{2\pi}\,\x'(t)\times\vv(t)\cdot\vv'_P(t).$$\rm

\medskip

\noindent\it Proof\rm.\ \ To begin, we
apply Taylor's
theorem to $K$ and write:
$$\x(s)=\x(t)+(s-t)\T(t)+\frac{(s-t)^2}{2}\x''(t)+(s-t)^3\A_1(s,t),$$
where $\T(t)=\x'(t)$ is the (unit) tangent vector to $K$ at $\x(t)$
and $\A_1(s,t)$ is a smooth, bounded (independent of $\eps$)
vector-valued
function of $s$ and $t$.

Because the link and writhe integrals (even the partial ones)
are
invariant under shifting the intervals of integration
(i.e., adding different constants mod $L$ to $s$ and $t$),
we may, without loss of generality, assume that
$t=0$. Then we
can write:
$$\x(s)=\x + s\T+\frac{s^2}{2}\x''+s^3\A_1(s)$$
where $\x=\x(0)$, where $\T=\x'(0)$ is the unit tangent vector to $K$ at $\x$,
where $\x''=\x''(0)$, and where
$\A_1(s)$ is a smooth, bounded (independent of $\eps$ and uniformly in $t$)
vector-valued
function of $s$. As in Proposition 5.1, we will write $\aes$ for $\alpha(\x(s),\ye(0))$ in
what follows.

Similarly, we can write
$$\frac{d\x}{ds}=\T+s\x''+s^2\A_2(s),$$
and we recall that the other edge of the ribbon and its derivative are
given by
$$\ye(0)=\cos\eps\,\x+\sin\eps\,\vv$$
and
$$\left.\frac{d\ye}{dt}\right|_{t=0}=\cos\eps\,\T+\sin\eps\,\vv',$$
where $\vv=\vv(0)$ and $\vv'=\vv'(0)$.
Because we are differentiating $\vv$ as though it were a vector field in
$\reals^4$, the derivative here coincides with the covariant derivative
on $S^3$, rather than the left-invariant derivative
of section 4. Here and for the remainder of this section, until the statement of the
theorem, we will omit the subscript in the notation $\vv'_P$.

Using the notation of
Proposition 5.1,
we can express
$$F_\eps(s,0)= 
-\frac{\ph'(\aes)}{\sin\aes}
\left\Vert\ye(0),\frac{d\ye}{dt}(0),\x(s),\frac{d\x}{ds}\right\Vert$$
as $-{\ph'(\aes)/\sin\aes}$ times the determinant
$$\left\Vert\cos\eps\,\x+\sin\eps\,\vv,\ \cos\eps\,\T+
\sin\eps\,\vv',\ \x+s\T+\frac{s^2}{2}
\x''+s^3\A_1,\ \T+
s\x''+s^2\A_2\right\Vert.$$

\vfill
\eject
We
proceed to analyze the factor $\ph'(\aes)/\sin\aes$ in
front of the determinant, and the following four terms,
into which the determinant can be expanded:
$$\eqalign{
I&=\left\Vert\cos\eps\,\x,\ \cos\eps\,\T,\ \frac{s^2}{2}\x''+s^3\A_1,\ s\x''+s^2
\A_2\right\Vert\cr
II&=\left\Vert\cos\eps\,\x,\ \sin\eps\,\vv',\ s\T+\frac{s^2}{2}
\x''+s^3\A_1,\ \T+s\x''+s^2\A_2\right\Vert\cr
III&=\left\Vert\sin\eps\,\vv,\ \cos\eps\,\T,\ \x+\frac{s^2}{2}
\x''+s^3\A_1,\ s\x''+s^2\A_2\right\Vert\cr
IV&=\left\Vert\sin\eps\,\vv,\ \sin\eps\,\vv',\ \x+s\T+\frac{s^2}{2}\x''+
s^3\A_1,\ \T+s\x''+s^2\A_2\right\Vert.}$$

First,
we derive an expansion of
$\ph'(\aes)/\sin\aes$ in powers of $s$ and $\eps$.
To begin, recall that, since $\x(s)$ is a curve on $S^3$ and is parametrized by
arclength, we have $\ip{\x}{\x}=1$, $\ip{\x}{\T}=0$, $\ip{\T}{\T}=1$ and 
$\ip{\x}{\x''}=-1$ (the last equation comes from differentiating $\ip{\x}{\T}=0$). 
Using these observations, we derive
$$\eqalign{
\cos\aes&=\ip{\x(s)}{\ye(0)}\cr
&=\ip{\x+s\T+\frac{s^2}{2}\x''+s^3\A_1\,}{\ \cos\eps\,\x+\sin\eps\,\vv}\cr
&=\cos\eps-\frac{s^2}{2}\cos\eps+\frac{s^2}{2}\sin\eps\,\ip{\x''}{\vv}+s^3p_0\cr
&=1-\frac{\eps^2+s^2}{2}+s^3p_1+s^2\eps p_2+s\eps^2 p_3+\eps^3 p_4}$$
where, as before, $p_i$ stands for a function of $s$ and $\eps$
that is bounded for all $s$ and $\eps$, and smooth except perhaps for $s=\eps=0$.

Since $\sin^2\alpha=1-\cos^2\alpha$, we can conclude that
$$\eqalign{
\sin^2\aes&=\eps^2+s^2+s^3p_5+s^2\eps p_6+s\eps^2p_7+\eps^3p_8\cr
&=(\eps^2+s^2)\left(1+\frac{s^3}{\eps^2+s^2}p_5+\frac{s^2\eps}{\eps^2+s^2}p_6
+\frac{s\eps^2}{\eps^2+s^2}p_7+\frac{\eps^3}{\eps^2+s^2}p_8\right).}$$
Using the Taylor series $\sqrt{1+z}=1+\frac{1}{2}z+\cdots$ we can conclude that
$$\sin\aes=(\eps^2+s^2)^{1/2}\left(1+\frac{s^3}{\eps^2+s^2}p_9+
\frac{s^2\eps}{\eps^2+s^2}p_{10}
+\frac{s\eps^2}{\eps^2+s^2}p_{11}+\frac{\eps^3}{\eps^2+s^2}p_{12}\right).$$
Using the Taylor series \ $\arcsin z = z+\cdots$ we can conclude that
$$\alpha=(\eps^2+s^2)^{1/2}\left(1+\frac{s^3}{\eps^2+s^2}p_{13}+
\frac{s^2\eps}{\eps^2+s^2}p_{14}
+\frac{s\eps^2}{\eps^2+s^2}p_{15}+\frac{\eps^3}{\eps^2+s^2}p_{16}\right).$$
We combine this with the expansion of $\ph'(\alpha)=-1/(4\pi\alpha^2)+$ something
bounded so that
$$\ph'(\aes)=-\frac{1}{4\pi(\eps^2+s^2)}+p_{17},$$
and finally conclude that
$$\eqalign{
\frac{\ph'(\aes)}{\sin\aes}&=\left(\frac{-1}{4\pi(\eps^2+s^2)}+p_{17}\right)
\left(\frac{1}{(\eps^2+s^2)^{1/2}}\right)(1+sp_{18}+\eps p_{19})\cr
&=\frac{-1}{4\pi(\eps^2+s^2)^{3/2}}(1+sp_{20}+\eps p_{21}).}$$

The utility of this expression for $\ph'(\aes)/\sin\aes$ will become
apparent
when we multiply it by the determinants, integrate from $-\eps^p$ to
$\eps^p$,
and then take the limit as $\eps\to 0$. Because $(\eps^2+s^2)^{1/2}$ is larger than
either $s$ or $\eps$, we can see that whenever $a+b\ge 3$, the product of
$\ph'(\aes)/\sin\aes$ with
$s^a\eps^b$ will integrate to something comparable to $\eps^p$, and the
integral  
will go to zero as $\eps$ does.

Next, we will use the observation about $(\eps^2+s^2)^{1/2}$ from the preceding
paragraph to deal with the four determinants. The first
one,
$$I=\left\Vert\cos\eps\,\x,\ \cos\eps\,\T,\ \frac{s^2}{2}\x''+s^3\A_1,\ s\kappa\N+s^2\A_2\right\Vert,$$
clearly has a factor of $s^3$, so it will not contribute to our limit. Similarly,
the second one,
$$II=\left\Vert\cos\eps\,\x,\ \sin\eps\,\vv',\ s\T+\frac{s^2}{2}
\x''+s^3\A_1,\ \T+s\kappa\N+s^2\A_2\right\Vert,$$
has a factor of $s^2\sin\eps$ (since you can't use the $\T$ from both the third and
fourth rows), and so $II$ doesn't contribute to our limit, either.

Using the expansion $\cos\eps=1-\eps^2/2+\cdots$, we can express the third determinant,
$$III=\left\Vert\sin\eps\,\vv,\ \cos\eps\,\T,\ \x+\frac{s^2}{2}\x''+s^3\A_1,\ s\x''+s^2\A_2\right\Vert,$$
as the sum of two terms:
$$s\sin\eps\left\Vert\vv,\T,\x,\x''\right\Vert+s^2\sin\eps\ p_{22},$$
from which only the first term could contribute to our limit.

\vfill
\eject
Finally, the fourth determinant,
$$IV=\left\Vert\sin\eps\,\vv,\ \sin\eps\,\vv',\ \x+s\T+\frac{s^2}{2}\x''+
s^3\A_1,\ \T+s\x''+s^2\A_2\right\Vert,$$
can be decomposed as
$$\sin^2\eps\left\Vert\vv,\vv',\x,\T\right\Vert+s\sin^2\eps\ p_{23},$$
from which only the first term could contribute to our limit.

From our analysis so far, we conclude that
$$F_\eps(s,0) = \frac{1}{4\pi(\eps^2+s^2)^{3/2}}
(s\eps\left\Vert\vv,\T,\x,\x''\right\Vert
+\eps^2\left\Vert\vv,\vv',\x,\T\right\Vert
+Z(\eps,s)),$$ where
$$\lim_{\eps\to 0}\int_{-\eps^p}^{\eps^p}
\frac{Z(\eps,s)}{(\eps^2+s^2)^{3/2}}\,ds=0.$$

We are now
ready to calculate the limit of the integral:
$$\lim_{\eps\to 0}\int_{-\eps^p}^{\eps^p}
F_\eps(s,0)\,ds.$$
From the formula for the integrand given above,
this limit will equal
$$\lim_{\eps\to
0}\left(\left\Vert\vv,\T,\x,\x''\right\Vert\int_{-\eps^p}^{\eps^p}
\frac{s\eps}{4\pi(\eps^2+s^2)^{3/2}}\,ds +
\left\Vert\vv,\vv',\x,\T\right\Vert\int_{-\eps^p}^{\eps^p}
\frac{\eps^2}{4\pi(\eps^2+s^2)^{3/2}}\,ds\right).$$

The integrand in the first of these integrals is odd, so the integral is
always zero (and hence the limit of that term is zero). For the second term, we
will need the fact that (for $p<1$)
$$\lim_{\eps\to 0}\int_{-\eps^p}^{\eps^p}\frac{\eps^2}{(\eps^2+s^2)^{3/2}}\,ds = 2,$$
which one calculates using the substitution $x=s/\eps$ and the fact that the
anti-derivative of $1/(1+x^2)^{3/2}$ is $x/\sqrt{1+x^2}$.

We have thus reached our final conclusion, namely that
$$\lim_{\eps\to 0}\int_{t-\eps^p}^{t+\eps^p} F_\eps(s,t)
\,ds =
\frac{1}{2\pi}\left\Vert\vv,\vv',\x,\T\right\Vert
=\frac{1}{2\pi}(\T\times\vv\cdot\vv').$$
This completes the proof of Proposition 6.4.

\vfill
\eject

We can use
Propositions 6.3 and 6.4 and a little arithmetic to start
from the definitions of $H_{Lk}(t;\eps)$ and $H_{Wr}(t)$ and deduce:

\medskip

\noindent\bf Proposition 6.5\rm:
$$\lim_{\eps\to 0}H_{Lk}(t;\eps) = H_{Wr}(t) +
\frac{1}{2\pi}(\T\times\vv\cdot\vv').$$

\medskip

We integrate the expression in Proposition 6.5
with respect to $t$ from $0$ to $L$
to reach our final conclusion:

\medskip

\noindent\bf Theorem 6.6\rm:\
$$\Lk(\x,\y)=\Tw_P(\x,\vv)+ \Wr_P(\x).$$
\it  
In other words, \sc Link = Twist + Writhe\rm.

\medskip

\noindent\bf Example\rm.\ \ The simplest example of two
linked curves
on $S^3$ is a pair of great circles from the same Hopf fibration. We verify
Theorem 6.6 in this case. The
curve\\
 $\x(s)=[\cos s,\sin s,0,0]$ is a great circle parametrized by
arclength
as $s$ runs from $0$ to $L=2\pi$.
We
will take $\x$ as one edge of our ribbon.

%%% FIGURE 3 %%%
\begin{figure}[h!]
\center{\includegraphics[height=190pt]{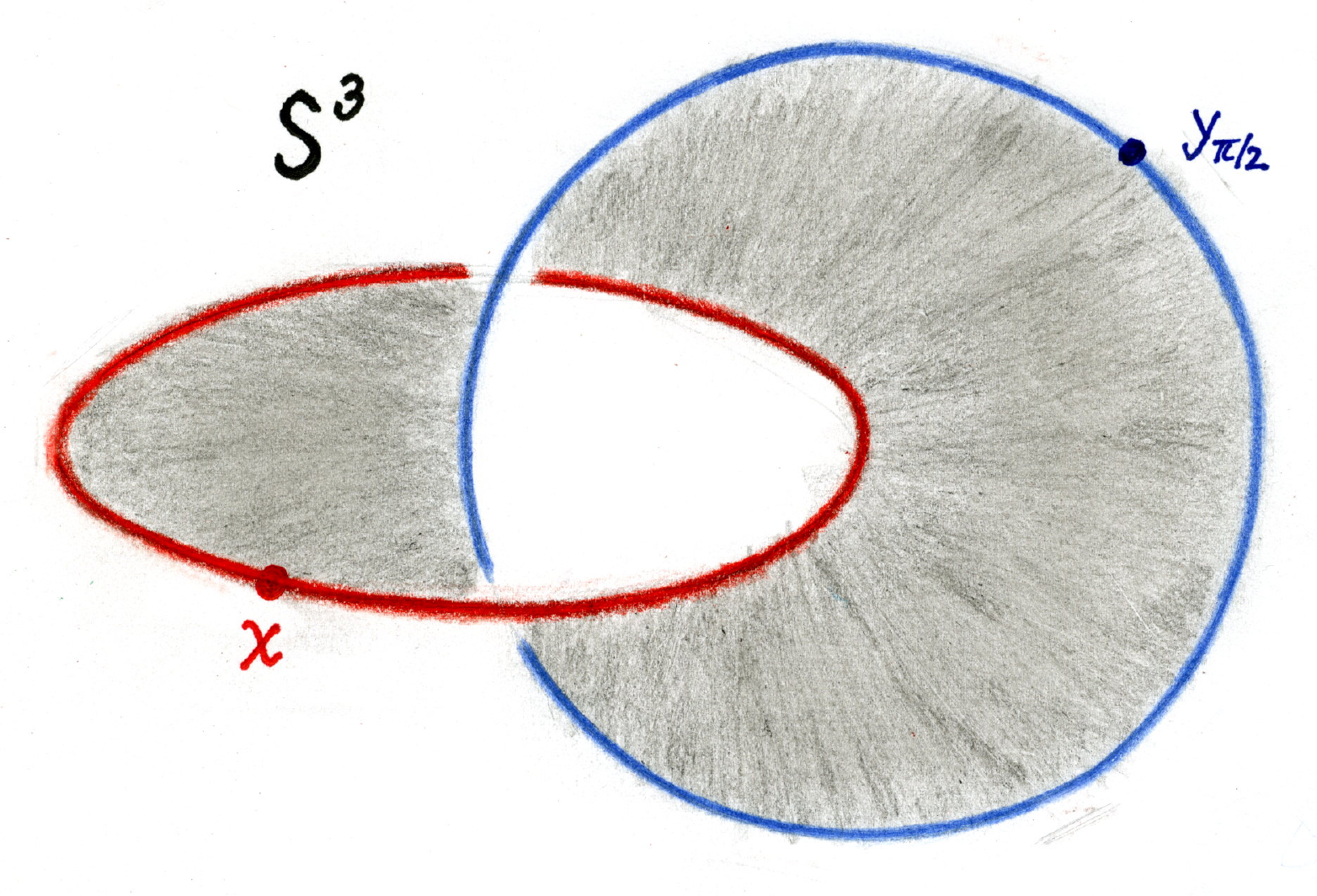}
\caption{The ribbon in the example}}
\end{figure}
%%%%%%%%%%%%%%%%%%%%

Let $\vv(t)=[0,0,\cos t,\sin t]$. Then $\vv$ is the restriction of a
left-invariant vector field to the great circle, and we will take the other
edge of our ribbon to be
$$\eqalign{
\ye(t)&=\cos\eps\,\x(t)+\sin\eps\,\vv(t)\cr
&= [\cos\eps\,\cos t\,,\,\cos\eps\,\sin t\,,\, \sin\eps\,\cos t\,,\,
\sin\eps\,\sin t]. }$$

If $\eps=\pi/2$, then $\ye(t)=[0,0,\cos t,\sin t]$, which is
the ``orthogonal'' great circle to $\x$ and we compute the linking
number of these two circles
as follows. Since $\ip{\x(s)}{\ye(t)}=0$ for all $s$ and $t$, we have
that $\alpha(\x(s),\y(t))=\pi/2$ for all $s$ and $t$.
Therefore, the linking integrand is given by
$$\eqalign{
\frac{d\x}{ds}\cdot P_{\x\ye}\left(\frac{d\ye}{dt}\times
\nabla_{\ye} \ph\right)&=\frac{-\ph'(\alpha)}{\sin\alpha}\,
\left\Vert\ye(t)\,,\,\frac{d\ye}{dt}\,,\,\x(s)\,,\,\frac{d\x}{ds}\right\Vert\cr
&=\frac{-\ph'(\pi/2)}{\sin(\pi/2)}\det\left|\begin{array}{cccc}
0 & 0 & \cos t & \sin t \\
0 & 0 & -\sin t& \cos t \\
\cos s& \sin s & 0 & 0\\
-\sin s& \cos s & 0 & 0 \end{array}\right|\cr
&=\frac{1}{4\pi^2}}$$
The integration takes place for $(s,t)\in[0,2\pi]\times [0,2\pi]$, so the formula for the
linking number of $\x$ and $\ye$ yields 1,
as expected.

To calculate the twist of our ribbon,
we note that $\T(s)=\x'(s)=[-\sin s,\cos s,0,0]$, and
$\vv'(s)=[0,0,-\sin s,\cos s]$
It is then easy to calculate that $\T\times\vv\cdot\vv'=[\x,\T,\vv,\vv']=1$
for all
$s$,
which gives us that the twist of the ribbon is
$$\Tw(\x,\vv)=\frac{1}{2\pi}\int_0^L\T(t)\times\vv(t)\cdot\vv'(t)\,dt
=\frac{1}{2\pi}\int_0^{2\pi} 1\,dt = 1.$$
To calculate the writhe of $\x$, we use the fact that $\x$ is
a geodesic, and so we have $P_{\x(s)\x(t)}\T(t)=\T(s)$. From this it is easy
to conclude that
$$\frac{d\x}{ds}\cdot
P_{\x(s)\x(t)}\left(\frac{d\x}{dt}\times
\nabla_{\x(t)}\ph\bigl(\alpha(\x(s),\x(t))\bigr)\right)=0$$
for all $s$ and $t$.
Therefore $$\Wr(\x)=0.$$
Theorem 6.6 then reads
$$\Lk(\x,\ye)=\Tw(\x,\vv)+\Wr(\x)= 1+0=1$$
as it should.

\vfill
\eject

\noindent\bf 7. Proof of \sc link = twist + writhe \bf in $H^3$\rm

The proof of \sc link = twist + writhe  in $H^3$ \rm is essentially a
repetition of
the parallel transport format proof in $S^3$, except for various
changes of sign and replacing trigonometric functions with
their corresponding hyperbolic ones. In this section, we highlight the
places where differences occur.

As in the first paper in this series, we view $H^3\subset\reals^{1,3}$,
the four-dimensional Minkowski space endowed with the inner product
$$\ip{\x}{\y}=x_0y_0-x_1y_1-x_2y_2-x_3y_3$$
so that
$$H^3=\{\x\in\reals^4\,|\,\ip{\x}{\x}=1\ \mbox{\rm and}\ x_0>0\}.$$
We
reserve the notation $\vv\cdot \w$ for the induced
inner product on $H^3$, namely for\\ $\vv,\w\in T_\x H^3$, we define
$\vv\cdot \w=-\ip{\vv}{\w}$.
Because the tangent vectors are spacelike, this inner product provides
$H^3$ with a Riemannian
metric which is complete and has constant curvature $-1$.

If $\x\in H^3$, and $\u,\vv\in T_\x H^3$, then we have
$$\u\times\vv = \det\left\vert \begin{array}{rccc}
x_0 & x_1 & x_2 & x_3 \\
u_0 & u_1 & u_2 & u_3 \\
v_0 & v_1 & v_2 & v_3 \\
-\xh_0 & \xh_1 & \xh_2 & \xh_3 \end{array} \right\vert .$$
Then for $\w\in T_\x H^3$,
the triple product $\u\times\vv\cdot\w = \Vert\x,\u,\vv,\w\Vert$.

For geodesics and the distance function, we will have
$$G(t)=\cosh t\ \x + \sinh t\ \vv$$
for the unit-speed geodesic through $\x$ in the direction of $\vv\in
T_\x H^3$, and the geodesic distance between $\x$ and $\y$ in $H^3$ will
satisfy $\cosh\alpha=\ip{\x}{\y}$. We have
$$\nabla_\y\alpha(\x,\y)=\frac{\cosh\alpha\ \y-\x}{\sinh\alpha}.$$
Except for the change in the inner product, the formula for parallel
transport remains the same: the result of parallel transport in $H^3$
of $\vv$ from $\y$ to $\x$ is
$$P_{\x\y}(\vv) =\vv - \frac{\ip{\x}{\vv}}{1+\ip{\x}{\y}}(\x+\y).$$

\vfill
\eject
Armed with these changes, and with the appropriate choice of
$\ph(\alpha)=\csch(\alpha)/(4\pi)$, the proofs of Lemma 6.1, Corollary 6.2,
Proposition 6.3 (where the biggest change is to have $\sinh \alpha_\eps$
rather than $\sin \alpha_\eps$ in the denominator) and Proposition 6.4
proceed in the hyperbolic space case essentially without change from
the spherical case.

We are then led to the conclusion of Proposition 6.5,
$$\lim_{\eps\to 0}H_{Lk}(t;\eps) = H_{Wr}(t) +
\frac{1}{2\pi}(\T\times\vv\cdot\vv').$$
And once again, we  
define the writhe of the $\x$ edge of our
ribbon as
$$Wr(\x)=\int_0^L\!\!\!\int_0^L \frac{d\x}{ds}\cdot
P_{\x(s)\x(t)}\left(\frac{d\x}{dt}\times
\nabla_{\x(t)}\ph\bigl(\alpha(\x(s),\x(t))\bigr)\right)\,ds\,dt,$$
and the twist of our ribbon as
$$Tw(\x,\vv)=\frac{1}{2\pi}\int_0^L (\T\times\vv\cdot\vv')\,dt.$$

Finally, we integrate the expressions from the hyperbolic
version of Proposition 6.5
with respect to $t$ from $0$ to $L$
to reach our final conclusion:

\medskip

\noindent\bf Theorem 7.1\rm:\
$$Lk(\x,\ye)=Tw(\x,\vv)+ Wr(\x).$$
\it  
In other words, \sc Link = Twist + Writhe.\rm

\medskip

\noindent\bf Example\rm.\ \ A simple example of a ribbon
in $H^3\subset\reals^{1,3}$ has as one edge the circle
$$\x(s)=[\sqrt{2},\cos s,\sin s, 0]$$
in $\reals^{1,3}$. The unit tangent vector to this curve is 
$$\T=[0,-\sin s,\cos s,0],$$
and we can choose the vector field
$$\vv(s)=\frac{1}{\sqrt{1-\frac{1}{2}\cos^2 s}}\left[\frac{\cos
s}{\sqrt{2}},\cos^2 s,\cos s\sin s,\sin s\right]$$
along $\x$. Clearly, $\ip{\x'(s)}{\vv(s)}=\ip{\T(s)}{\vv(s)}=0$ for all $s$,
and $\ip{\vv(s)}{\vv(s)}=-1$, so $\vv$ is a unit vector perpendicular to $T$.
We can make the ribbon by choosing the other edge to be the curve given
by $\ye(s) = \cosh\eps\ \x(s) + \sinh\eps\ \vv(s)$.

%%% FIGURE 3 %%%
\begin{figure}[h!]
\center{\includegraphics[height=160pt]{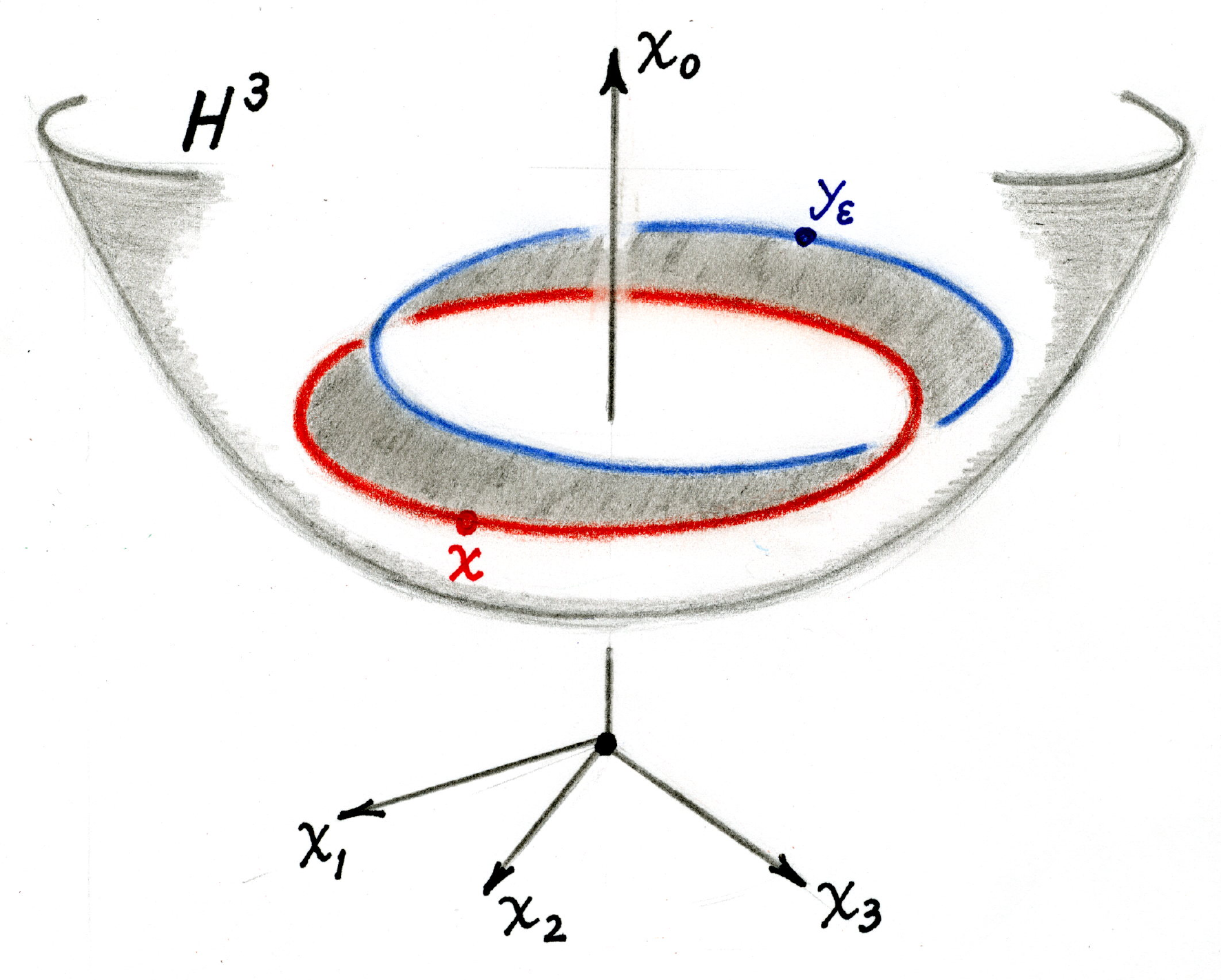}
\caption{The ribbon in the example}}
\end{figure}
%%%%%%%%%%%%%%%%%%%%

By looking at the projections of the $\x$ and $\ye$ curves into $\reals^3$
(ignoring the first coordinates), it's easy to see that these curves have
linking number $-1$. The writhing integrand of the $\x$ curve is easily
seen to be zero, since the writhing
integrand
is given by
$$-\frac{\ph'(\alpha)}{\sinh\alpha}\left\Vert\x(t),\frac{d\x}{dt},\x(s),
\frac{d\x}{ds}\right\Vert,$$ and the determinant is zero because the last component of
each vector in the determinant is zero.

For the twist of the ribbon, we must calculate $\T\times\vv\cdot\vv'$, which
is given by the determinant
$$\left\Vert\x(s),\vv(s)\,,\,\frac{d\vv}{ds},\T(s)\right\Vert,$$ and calculating this  
determinant yields
$$\T\times\vv\cdot\vv'=\frac{\sqrt{2}}{\cos^2 s-2}$$
So we can calcuate that
$$\eqalign{
\Tw(\x,\vv)&=\frac{1}{2\pi}\int_0^{2\pi} (\T\times\vv\cdot\vv')\,ds\cr
&=\frac{1}{2\pi}\int_0^{2\pi}\frac{\sqrt{2}}{\cos^2 s-2}\,ds=-1.}$$
Theorem 7.1 then reads
$$\Lk(\x,\ye)=\Tw(\x,\vv)+\Wr(\x)=-1+0=-1$$
as it should.

\bigskip

\noindent\bf 8. Helicity of vector fields on $S^3$ and $H^3$\rm

Lodewijk Woltjer introduced in 1958 the notion of ``helicity'' of a vector
field $\vv$ defined on a domain $\Omega$ in Euclidean 3-space,
$$\He(\vv)=\int_{\Omega\times\Omega}\vv(\x)\times\vv(\y)\cdot
\frac{\x-\y}{4\pi|\x-\y|^3}\,d\x\,
d\y,\eqno(8.1)$$
as an invariant during ideal magnetohydrodynamic evolution of plasma fields. Keith
Moffatt (1969), recognizing that this quantity measures the extent to which the field lines
of $\vv$ wrap and coil around one another, named it ``helicity'' and showed that
Woltjer's original formula could be written in the above form.

If $\vv$ is a smooth vector field on $\reals^3$ with compact support, then the above formula for
its helicity can be written succinctly as
$$\He(\vv)=\int_{\reals^3}\BS(\vv)(\y)\cdot\vv(\y)\,d\y,\eqno(8.2)$$
where we recall that $\BS(\vv)$ denotes the magnetic field due to the steady current flow 
$\vv$.

This is how Woltjer originally presented his invariant, $\displaystyle{\int \A\cdot \B\,d\x}$, with the
role of $\vv$ played by the magnetic field $\B$ and the role of $\BS(\vv)$ played by its
vector potential $\A$.

We use (8.2) to define the helicity of a vector field $\vv$ on $S^3$ or $H^3$,
and then immediately obtain explicit integral formulas from Theorem 2.1.

\medskip

\noindent\bf Theorem 8.3\rm. \sc Helicity integrals in $S^3$ and
$H^3$\it.

\rm(1) \it On $S^3$, in left-translation format:
$$\eqalign{
\He(\vv)=\int_{S^3\times
S^3}L_{\y\x^{-1}}\vv(\x)&\times\vv(\y)\cdot\nabla_\y
\ph(\x,\y)\,d\x\,d\y -\frac{1}{4\pi^2}\int_{S^3\times S^3} L_{\y\x^{-1}}\vv(\x)
\cdot\vv(\y)\,d\x\,d\y\cr
&+2\int_{S^3\times S^3}\nabla_\y (L_{\y\x^{-1}}\vv(\x)\cdot\nabla_\y
\ph_1(\x,\y))\cdot\vv(\y)\,d\x\,d\y,}$$
where $\ph(\alpha)=(\pi-\alpha)\cot\alpha/(4\pi^2)$ and
$\ph_1(\alpha)=-\alpha(2\pi-\alpha)/(16\pi^2)$.

\rm(2) \it On $S^3$ in parallel transport format:
$$\He(\vv)=\int_{S^3\times S^3}P_{\y\x}\vv(\x)\times
\vv(\y)\cdot\nabla_\y\ph(\x,\y)\,d\x\,d\y,$$
where $\ph(\alpha)=(\pi-\alpha)\csc\alpha/(4\pi^2)$.

\rm(3) \it On $H^3$ in parallel transport format:
$$\He(\vv)=\int_{H^3\times H^3}P_{\y\x}\vv(\x)\times
\vv(\y)\cdot\nabla_\y\ph(\x,\y)\,d\x\,d\y$$
where $\ph(\alpha)=\csch\alpha/(4\pi)$.\rm

\medskip

In formula (1), if $\vv$ is divergence-free, then the third integral in the definition of
$\He(\vv)$ vanishes, and this formula then resembles the linking formula (1) of Theorem 1.1.
Formulas (2) and (3) already resemble the corresponding linking formulas of Theorem~1.1.

In formulas (1) and (2), if the smooth vector field $\vv$ on $S^3$ is divergence-free,
then its helicity is the same as its \it asymptotic \rm(or \it mean\rm) \it Hopf invariant\rm,
as defined by Arnold (1974), and is invariant under the group of volume-preserving diffeomorphisms
of $S^3$.

In formula (3), we assume that $\vv$ has compact support in order to guarantee convergence of the
integral.

\bigskip

\noindent\bf 9. Upper bounds for helicity in $\reals^3$, $S^3$ and $H^3$\rm.

Let $\Omega$ be a compact, smoothly bounded subdomain of $\reals^3$, $S^3$ or $H^3$, and let
$\vv$ be a smooth vector field defined on $\Omega$. 
Thinking of $\vv$ as a current flow, its magnetic field
$\BS(\vv)$ is defined by the same formulas as in Theorem 2.1, except that the integration is carried out
only over $\Omega$.

For uniformity of approach, we ignore the left-translation format on $S^3$ and write
$$\BS(\vv)(\y)=\int_\Omega P_{\y\x}\vv(\x)\times
\nabla_\y\ph_0(\x,\y)\,d\x,\eqno(9.1)$$
where 
$$\begin{array}{lll}
\mbox{\rm in $\reals^3$} & \mbox{\rm we have}\ \ph_0(\alpha)=-\displaystyle{\frac{1}{4\pi\alpha}}
&\mbox{\rm so $\Delta\ph_0=\delta$;}\\ \\
\mbox{\rm in $S^3$} & \mbox{\rm we have}\ \ph_0(\alpha)=\displaystyle{-\frac{1}{4\pi^2}}
(\pi-\alpha)\csc\alpha &\mbox{\rm so $\Delta\ph_0-\ph_0=\delta$;}\\ \\ 
\mbox{\rm in $H^3$} & \mbox{\rm we have}\ \ph_0(\alpha)=\displaystyle{-\frac{1}{4\pi}}\csch\alpha
&\mbox{\rm so $\Delta\ph_0+\ph_0=\delta$.}
\end{array}
$$

The magnetic field $\BS(\vv)$ is defined throughout the ambient space. It is continuous
everywhere, but its first derivatives suffer a discontinuity as one crosses the boundary
of $\Omega$. This is a familiar situation from electrodynamics in Euclidean 3-space.

In what follows, we will restrict $\BS(\vv)$ to $\Omega$, and ignore its behavior outside
this domain.

Let $\VF(\Omega)$ denote the space of all smooth vector fields on $\Omega$, with the $L^2$
inner product
$$\ip{\vv}{\w}=\int_\Omega \vv\cdot\w\,d\vol,$$
and associated \it energy \rm $\ip{\vv}{\vv}$ and \it norm \rm $|\vv|=\ip{\vv}{\vv}^{1/2}$.

We seek a bound for the energy or norm of the output magnetic field $\BS(\vv)$ on
$\Omega$ in terms of the input current flow $\vv$. Or to put it another way, we seek an upper
bound for the $L^2$-operator norm of the Biot-Savart operator,
$$\BS\colon\VF(\Omega)\to\VF(\Omega),$$
in terms of the geometry of the underlying domain $\Omega$.

As a consequence, we will determine an upper bound for the helicity
$\He(\vv)=\ip{\BS(\vv)}{\vv}$
of the vector field $\vv$ in terms of its energy $\ip{\vv}{\vv}$ and the geometry of $\Omega$.

\medskip

\noindent\bf Theorem 9.2\it.\ \ Let $\Omega$ be a compact, smoothly
bounded subdomain of $\reals^3$, $S^3$ or $H^3$ and let $R=R(\Omega)$ be
the radius of a ball in that space having the same volume as $\Omega$.
Let $\vv$ be a smooth vector field defined on $\Omega$. Then
$$|\BS(\vv)|\le N(R)|\vv|,$$
where
$$\begin{array}{ll}
\mbox{\it in $\reals^3$} & \mbox{\it we have}\ N(R)=R\\ \\
\mbox{\it in $S^3$} & \mbox{\it we have}\ N(R)=\displaystyle{\frac{1}{\pi}(2(1-\cos R)+(\pi-R)\sin R)}\\ \\
\mbox{\it in $H^3$} & \mbox{\it we have}\ N(R)=\sinh R.
\end{array}
$$

%%% FIGURE  %%%
\begin{figure}[h!]
\center{\includegraphics[height=160pt]{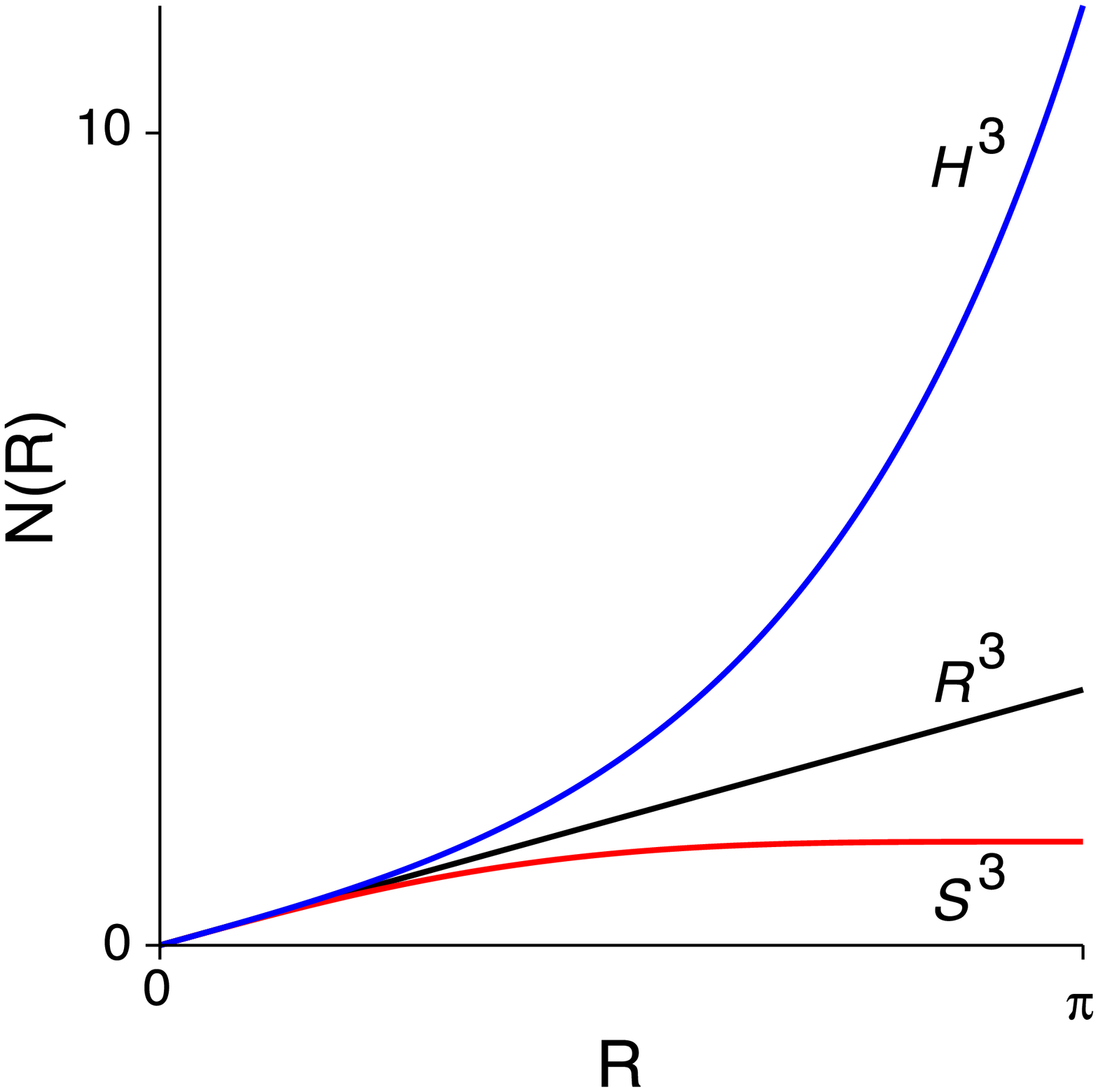}
\caption{$N(R)$ for $H^3$, $\reals^3$ and $S^3$}}
\end{figure}
%%%%%%%%%%%%%%%%%%%%

\noindent It follows immediately that the helicity $\He(\vv)=\ip{\BS(\vv)}{\vv}$ is
bounded by
$$|\He(\vv)|\le N(R)|\vv|^2.$$\rm

\medskip

In $\reals^3$, the overestimate $N(R)=R$ for the norm of the
Biot-Savart operator grows like the \it cube root \rm of the volume
$\frac{4}{3}\pi R^3$ of $\Omega$.

By contrast, in $H^3$ the overestimate $N(R)=\sinh R$ for the norm of
the Biot-Savart operator grows like the \it square root \rm of the volume
$2\pi(\sinh R\cosh R-R)$ of $\Omega$.

\medskip

\noindent\bf Setting up for the proof of Theorem 9.2.\rm

To begin, let $\psi(\alpha)$ be a real-valued function of the real variable
$\alpha>0$, where we think of $\alpha$ as distance from a fixed point (and
on $S^3$, we have the additional condition $0<\alpha\le\pi$). 
Assume $\psi$ has the property that
$$N_\Omega(\psi):= \max_\y \int_{\x\in\Omega}|\psi(\x,\y)|\,d\x$$
is finite, where as usual we write $\psi(\x,\y)$ as an abbreviation for
$\psi(\alpha(\x,\y))$. We note explicitly that the point $\y$ need \it not \rm
be chosen in $\Omega$.

\medskip

\noindent\bf Proposition 9.3\it.\ \ Under the above circumstances, the operator
$$T_\psi\colon \VF(\Omega)\to \VF(\Omega)$$
defined by 
$$T_\psi(\vv)(\y)=\int_\Omega P_{\y\x}\vv(\x)\times\psi(\x,\y)\nabla_\y\alpha(\x,\y)\,d\x$$
is a bounded operator with respect to the $L^2$-norm, and
$$|T_\psi(\vv)|\le N_\Omega(\psi)|\vv|.$$\rm

\medskip

The proof of this proposition in the $\reals^3$ case can be found in our (2001) paper,
Lemma 3 on pages 897 and 898. The argument there follows along the lines of the
usual Young's inequality proof that convolution operators on spaces of scalar-valued
functions are bounded; see Folland (1995) page 9, or Zimmer (1990) Proposition B.3 on
page 10. The proof carries over to the $S^3$ and $H^3$ cases with virtually no changes.

We want to apply Proposition 9.3 to the Biot-Savart operator (9.1), which we
write as 
$$
\BS(\vv)(\y) = \int_\Omega P_{\y\x}\vv(\x)\times\ph_0'(\x,\y)\nabla_\y\alpha(\x,\y)\,d\x
=T_{\ph_0'}(\vv)(\y),
\eqno(9.4)$$
where
$$\begin{array}{ll}
\mbox{\it in $\reals^3$} & \mbox{\it we have}\ \ph_0'(\alpha)=
\displaystyle\frac{1}{4\pi\alpha^2}\,; \\ \\
\mbox{\it in $S^3$} & \mbox{\it we have}\ \ph_0'(\alpha)=
\displaystyle{\frac{1}{4\pi^2}}(\csc\alpha+(\pi-\alpha)\csc\alpha\cot\alpha)\,; \\ \\
\mbox{\it in $H^3$} & \mbox{\it we have}\ \ph_0'(\alpha)=
\displaystyle{\frac{1}{4\pi}}\csch\alpha\coth\alpha.
\end{array}
$$

Then by Proposition 9.3 we have 

\medskip

\noindent\bf Proposition 9.5\rm.\ \ $|\BS(\vv)|\le N_\Omega(\ph_0')|\vv|$.

\medskip

We turn next to estimating $N_\Omega(\ph_0')$.

\medskip

\noindent\bf Lemma 9.6\it.\ \ If $\psi(\alpha)$ is a positive, decreasing function of $\alpha$, 
then 
$$N_\Omega(\psi)=\max_\y\int_{\x\in\Omega}\psi(\x,\y)\,d\x$$
is maximized over all subdomains $\Omega$ of $\reals^3$, $S^3$ or $H^3$ having
a given volume when $\Omega$ is a round ball
and $\y$ is its center.\rm

\medskip

We leave the proof of this, as well as that of the next elementary lemma, to the reader.

\medskip

\noindent\bf Lemma 9.7\it.\ \ The functions 
$$\begin{array}{ll}
\ph_0'(\alpha)=\displaystyle{\frac{1}{4\pi\alpha^2}}&\quad \alpha\in(0,\infty)\\ \\
\ph_0'(\alpha)= \displaystyle{\frac{1}{4\pi^2}}(\csc\alpha+(\pi-\alpha)\csc\alpha\cot\alpha) 
&\quad \alpha\in(0,\pi]\\ \\
\ph_0'(\alpha)=\displaystyle{\frac{1}{4\pi}}\csch\alpha\coth\alpha&\quad \alpha\in(0,\infty)
\end{array}$$

\noindent are decreasing functions of $\alpha$ on their respective domains.\rm

\medskip

In view of Lemmas 9.6 and 9.7, we next compute $N_\Omega(\ph_0')$, where $\Omega$ is a round ball
of radius $R$ in $\reals^3$, $S^3$ or $H^3$, and $\ph_0'(\alpha)$ is as given above.
We use the shorthand $N(R)=N_\Omega(\ph_0')$.

\medskip

\noindent\bf Proposition 9.8\it.\ \ Let $\Omega$ be a round ball of radius $R$ in 
$\reals^3$, $S^3$ or $H^3$. Then

$$\begin{array}{ll}
\mbox{\it in $\reals^3$ we have}\ &N(R)=R\\ \\
\mbox{\it in $S^3$ we have}\ &N(R)= \displaystyle{\frac{1}{\pi}}(2(1-\cos R)+(\pi-R)\sin R) \\ \\
\mbox{\it in $H^3$ we have}\ &N(R)=\sinh R.
\end{array}$$\rm

\medskip

\noindent\it Proof\rm.\ \ We give the proof in $S^3$ and leave the other two cases to 
the reader.
$$\eqalign{
N(R) &= \int_{\alpha=0}^R \ph_0'(\alpha)\,4\pi\sin^2\alpha\,d\alpha\cr
&=\int_{\alpha=0}^R \frac{1}{4\pi^2}(\csc\alpha+(\pi-\alpha)\csc\alpha\cot\alpha)\,4\pi\sin^2\alpha\,
d\alpha\cr
&= \frac{1}{\pi}\int_{\alpha=0}^R(\sin\alpha+(\pi-\alpha)\cos\alpha)\,d\alpha\cr
&= \frac{1}{\pi}(2(1-\cos R)+(\pi-R)\sin R).
}$$

\medskip

\noindent\bf Remark\rm.\ \ If we put $R=\pi$, then 
we get $N_{S^3}(\ph_0')=4/\pi$, in which case
$$|\BS(\vv)|\le N_{S^3}(\ph_0')|\vv|=\frac{4}{\pi}|\vv|,$$
for smooth vector fields $\vv$ defined on the entire 3-sphere.

We contrast this with the sharp estimate
$$|\BS(\vv)|\le\frac{1}{2}|\vv|,$$
with equality if and only if $\vv$ is a vector field of constant length tangent to a left
or right Hopf fibration of $S^3$. See our (2008) paper for details.

\vfill
\eject

\noindent\it Proof of Theorem \rm 9.2.\ \ By Lemma 9.7, the functions $\ph_0'(\alpha)$
in $\reals^3$, $S^3$ and $H^3$ are all positive and decreasing. Hence by Lemma 9.6, the quantity
$$N_\Omega(\ph_0')=\max_{\y}\int_{\x\in\Omega}\ph_0'(\x,\y)\,d\x$$
is maximized over subdomains $\Omega$ having a given volume when $\Omega$ is a round 
ball and $\y$ is its center. The values of $N_\Omega(\ph_0')$ in that case were 
calculated in Proposition 9.8, and inserting them into the estimate
$$|\BS(\vv)|\le N_\Omega(\ph_0')|\vv|$$
of Proposition 9.5, we get Theorem 9.2.

\bigskip

\noindent\bf 10. Hodge decomposition of vector fields\rm

In this section we collect, without proof, some information about the topology of compact subdomains
in $\reals^3$, $S^3$ and $H^3$, and about the structure of the space of vector fields on such domains.  
The reader will find the details in our (2002) paper.

Let $\Omega$ be a compact, smoothly bounded domain in  $\reals^3$, $S^3$ or $H^3$,  and 
$\VF(\Omega)$ the space of all smooth vector fields on $\Omega$,  with the $L^2$ inner product and
associated energy and norm, as defined in the preceding section.

Let $\K(\Omega)\subset\VF(\Omega)$ denote the subspace consisting of vector fields which are 
divergence-free and tangent to the boundary of $\Omega$,
$$\K(\Omega)  =  \{\vv\in\VF(\Omega) \colon \nabla\cdot\vv =0 ,\  \vv\cdot\n=  0\} ,$$
where $\n$ denotes the unit outward normal vector field along the boundary $\partial\Omega$
of  $\Omega$ .  These vector fields are just the incompressible fluid flows within a 
bounded domain, and in real life are naturally tangent to the boundary.  In the 
traditional passage from geometric knot theory to fluid dynamics, a knot is modeled 
by such a flow within a tubular neighborhood of itself, and the flows are then called 
\it fluid knots\rm,  accounting for the  ``K''  in the notation $\K(\Omega)$ .

Let  $\G(\Omega) \subset \VF(\Omega)$  denote the subspace of \it gradient fields\rm,
$$\G(\Omega)  =  \{\vv\in\VF(\Omega)\colon  \vv  = \nabla\ph\ \mbox{\rm for some smooth function}
\ \ph\colon\Omega\to\reals\}.$$ 
Then we have an $L^2$-orthogonal direct sum decomposition  
$$\VF(\Omega)=\K(\Omega)\oplus\G(\Omega).\eqno(10.1)$$
The spaces $\VF(\Omega)$, $\K(\Omega)$ and $\G(\Omega)$ are all infinite-dimensional.

Let $\HK(\Omega)\subset\K(\Omega)$ denote the subspace of vector fields which are not only
divergence-free and tangent to the boundary, but also curl-free,
$$\HK(\Omega) =\{\vv\in\VF(\Omega)\colon \nabla\cdot\vv=0,\ \nabla\times\vv=
\zero,\ \vv\cdot\n=0\}.$$
We call the elements of  $\HK(\Omega)$ \it harmonic knots\rm. The subspace 
$\HK(\Omega)$  is
finite-dimensional, and isomorphic to  $H_1(\Omega)$,  the one-dimensional homology
of $\Omega$  with real coefficients.  

The orthogonal decomposition (10.1), when further refined, yields the \it Hodge decomposition
\rm of $\VF(\Omega)$;  see our (2002) paper for details.

Let $\Omega^\dast$ denote the closure of the complement of $\Omega$ in $\reals^3$,
$S^3$ or $H^3$. Let $g$ denote the \it total genus \rm of $\partial\Omega$, that is, the sum of the genera
of its components. Then, using real coefficients, $H_1(\partial\Omega)$ is a $2g$-dimensional vector
space, while $H_1(\Omega)$ and $H_1(\Omega^\dast)$ are each $g$-dimensional, and we have the direct sum
decomposition
$$H_1(\partial\Omega)= \ker(H_1(\partial\Omega)\to H_1(\Omega)) + \ker(H_1(\partial\Omega)
\to H_1(\Omega^\dast)),\eqno(10.2)$$
where the above homomorphisms are induced by the inclusions $\partial\Omega\subset\Omega$ and
$\partial\Omega\subset\Omega^\dast$.

Let  $a_1,a_2,\ldots,a_g$ be a basis for  $\ker(H_1(\partial\Omega)\to H_1(\Omega))$,
and $b_1,b_2,\ldots,b_g$ a basis for $\ker(H_1(\partial\Omega)\to H_1(\Omega^\dast))$. 

If  $\vv\in\HK(\Omega)$,  then, since $\vv$  is curl-free,  its circulation  
$$\Circ(\vv,\gamma) = \int_\gamma(\vv(\x(t))\cdot\frac{d\x}{dt}\,dt$$
about any curve $\gamma$ in $\Omega$ depends only on the homology class of $\gamma$.  
So we can denote this circulation by $\Circ(\vv,[\gamma])$. 

With this notation, the real numbers $\Circ(\vv,a_1),\Circ(\vv,a_2),\ldots,\Circ(\vv,a_g)$
are all zero, since the homology classes $a_i$ on $\partial\Omega$ bound in $\Omega$.  By contrast,

\noindent (10.3) The real numbers $\Circ(\vv,b_1),\Circ(\vv,b_2),\ldots,\Circ(\vv,b_g)$ are in general
not zero, and in fact define an \it isomorphism \rm of $\HK(\Omega)\to\reals^g$.

\vfill
\eject
\noindent\bf 11. Spectral geometry of the curl operator in $\reals^3$, $S^3$ and
$H^3$\rm

As before let $\Omega$ be a compact, smoothly bounded subdomain of $\reals^3$, $S^3$
or $H^3$, and $\VF(\Omega)$ the infinite-dimensional space of smooth 
vector fields on $\Omega$ with the $L^2$ inner product. 

Now we are interested in curl eigenfields on $\Omega$, that is, vector
fields $\vv$ on $\Omega$ which satisfy $\nabla\times\vv = \lambda\vv$ for $\lambda\ne 0$.
In $\reals^3$, these fields are used to model stable plasma flows; see our (1999) paper.

Curl eigenfields exist for every value of $\lambda$. For example, in $\reals^3$, if
$$\vv=\sin\lambda z\i +\cos\lambda z\j,$$
then $\nabla\times\vv=\lambda\vv$.

We want to constrain the choice of vector fields $\vv$ by interior and boundary
conditions which guarantee that the curl operator on $\VF(\Omega)$ 
will have a discrete spectrum, 
while at the same time being reasonable
for physical applications. Then we want to find a lower bound for the absolute values of the
nonzero eigenvalues.

To begin, we will restrict our attention to the subspace $\K(\Omega)$ of fluid knots, 
discussed in section 10. The vector fields in $\K(\Omega)$ are divergence-free and 
tangent to the boundary of $\Omega$. Since a curl eigenfield $\vv$ with nonzero eigenvalue
$\lambda$ is automatically divergence-free, the only real constraint here is that of 
tangency to the boundary.

Let $\CK(\Omega)\subset K(\Omega)$ denote the subspace of vector fields whose curl lies
in $K(\Omega)$. Any eigenfield of the curl operator in $K(\Omega)$ must lie in $\CK(\Omega)$, 
so restricting our attention to $\CK(\Omega)$ is no further constraint.

\medskip

\noindent\bf Lemma 11.1\it.\ \ A vector field $\vv\in\K(\Omega)$ lies in the subspace 
$\CK(\Omega)$ if and
only if the circulation of $\vv$ around small
loops on $\partial\Omega$
vanishes.\rm

\medskip

\noindent\it Proof\rm.\ \ The circulation of $\vv$ around a small loop on $\partial\Omega$ 
equals the flux of $\nabla\times\vv$ through the small disk bounded by that loop. If this
flux is zero for all such loops, 
then the normal component of $\nabla\times\vv$ along $\partial\Omega$ must be zero, telling us that
$\nabla\times\vv$ is tangent to $\partial\Omega$, and hence that $\vv\in\CK(\Omega)$.

\noindent\bf Remarks\rm. 

(1) If the ciculation of $\vv$ vanishes around small loops on $\partial\Omega$, then it 
also vanishes around homologically trivial loops there.

(2) Any divergence-free vector field on $\Omega$ which 
vanishes on $\partial\Omega$ must lie in $\CK(\Omega)$.

\medskip

The kernel of the map \ $\curl\colon\CK(\Omega)\to\K(\Omega)$ consists of vector fields on $\Omega$
which are divergence-free, curl-free, and tangent to the boundary. These are the harmonic knots
$\HK(\Omega)$ introduced in section 10.

\medskip

Since we are interested in the spectral theory of the curl operator, 
we would like to know when $\curl\colon\CK(\Omega)\to\K(\Omega)$ is self-adjoint with respect to the $L^2$
inner product, that is, when can we promise that
$$\ip{\nabla\times\vv}{\w}=\ip{\vv}{\nabla\times\w}$$
for vector fields $\vv$ and $\w$ in $\CK(\Omega)$?

\medskip

\noindent\bf Lemma 11.2\it.\ \ Suppose that $\Omega$ is simply connected, or equivalently, that all 
the components of $\partial\Omega$ are 2-spheres. Then $\curl\colon\CK(\Omega)\to\K(\Omega)$ is self-adjoint.\rm

\medskip

\noindent\it Proof\rm.\ \ Recall the formula from vector calculus,
$$\nabla\cdot(\vv\times\w)= (\nabla\times\vv)\cdot\w-\vv\cdot(\nabla\times \w),$$
and integrate this over $\Omega$ to get
$$\eqalign{
\int_\Omega(\nabla\cdot(\vv\times\w)\,d\vol&=\int_\Omega (\nabla\times\vv)\cdot\w\,d\vol
-\int_\Omega \vv\cdot(\nabla\times \w)\,d\vol\cr
&=\ip{\nabla\times\vv}{\w}-\ip{\vv}{\nabla\times\w}.
}$$
The left-hand side equals 
$$\int_{\partial\Omega}(\vv\times\w)\cdot\n\,d\area,$$
and so the issue is to see when $\vv\times\w$ has zero flux through $\partial\Omega$.

By Lemma 11.1, we know that $\vv\in\CK$ if and only if $\vv$ has zero circulation
around small loops on $\partial\Omega$. Since $\Omega$ is simply connected,
$\partial\Omega$ is a union of 2-spheres, and so $\vv$ must have zero circulation 
around \it all \rm loops on $\partial\Omega$. 

But this means that the restriction of $\vv$ to $\partial\Omega$ is a gradient field on 
that surface. So we write $\vv|_{\partial\Omega}=\nabla f$, where $f\colon\partial\Omega
\to\reals$ is some smooth function, and the gradient is the ``surface gradient'' on 
$\partial\Omega$. 

Likewise, $\w|_{\partial\Omega}=\nabla g$ for some smooth function $g\colon\partial\Omega
\to\reals$. 

Now extend $f$ and $g$ to smooth functions $F$ and $G$ from $\Omega\to\reals$, and
consider the vector fields $\nabla F$ and $\nabla G$ defined on $\Omega$.

Because the cross product of two gradient fields is always divergence-free, that is, 
$$\nabla\cdot(\nabla F\times\nabla G)=(\nabla\times\nabla F)\cdot\nabla G
-\nabla F\cdot(\nabla\times\nabla G)=0,$$
we have
$$\int_{\partial\Omega}(\nabla F\times\nabla G)\cdot\n\,d\area=0.$$

Since $\nabla f$ and $\nabla g$ are, respectively, the tangential components of
$\nabla F$ and $\nabla G$ along $\partial\Omega$, we can write
$$\nabla F(\x)=\nabla f(\x)+ a(\x)\,\n(\x)\quad\mbox{\rm and}\quad
\nabla G(\x)=\nabla g(\x) + b(\x)\,\n(\x)$$
for $\x\in\partial\Omega$. From this we can see that
$$(\nabla F(\x)\times\nabla G(\x))\cdot\n(\x)=(\nabla f(\x)\times\nabla g(\x))\cdot
\n(\x)$$
along $\partial\Omega$, and hence
$$\eqalign{
\int_{\partial\Omega}(\vv\times\w)\cdot\n\,d\area 
&= \int_{\partial\Omega}(\nabla f\times\nabla g)\cdot\n\,d\area\cr
&= \int_{\partial\Omega}(\nabla F\times\nabla G)\cdot\n,d\area\cr
&=0.
}$$
Thus $\ip{\nabla\times\vv}{\w}=\ip{\vv}{\nabla\times\w}$ for all $\vv$ and $\w$ in
$\CK(\Omega)$, and so $\curl\colon\CK(\Omega)\to\K(\Omega)$ is self-adjoint
when $\Omega$ is simply connected, completing the proof of Lemma 11.2.

\medskip

When $\Omega$ is not simply connected, the operator 
$\curl\colon\CK(\Omega)\to\K(\Omega)$ is \it not \rm self-adjoint.
So we seek further sensible boundary conditions which will make 
this operator self-adjoint for any domain
$\Omega$.

To this end, let $\vv$ be a vector field in $\CK(\Omega)$. By definition of $\CK(\Omega)$,
the curculation of $\vv$ around all small loops on $\partial\Omega$ vanishes. But then the
circulation of $\vv$ around \it any \rm loop on $\partial\Omega$ depends only on the homology
class of that loop, giving us a linear map
$$\Circ(\vv)\colon H_1(\partial\Omega)\to\reals,$$
from the one-dimensional real homology of $\partial\Omega$ to the reals.

To see where this is heading, let $\Omega^\dast$ denote the closure of the complement of
$\Omega$ in $\reals^3$, $S^3$ or $H^3$, as in section 10, and let $g$ be the total
genus of $\partial\Omega$.

\vfill
\eject
Recall from section 10 the direct sum decomposition
$$H_1(\partial\Omega)=\ker(H_1(\partial\Omega)\to H_1(\Omega))+ \ker(H_1(\partial\Omega)
\to H_1(\Omega^\dast)),\eqno(10.2)$$
which splits a $2g$-dimensional space into two $g$-dimensional summands.

Now let $\AK(\Omega)\subset\CK(\Omega)\subset\K(\Omega)$ consist of all vector fields $\vv$ in $\CK(\Omega)$
whose circulation vanishes around any loop on $\partial\Omega$ which bounds in $\Omega^\dast$. 
The subspace $\AK(\Omega)$ has codimension $g$ in $\CK(\Omega)$.

From (10.3), we get
$$\AK(\Omega)\,\cap\,\HK(\Omega)=\{0\}.\eqno(11.3)$$
We call $\AK(\Omega)$ the space of \it Amp\`erian knots \rm because, by Amp\`ere's Law, the magnetic field
due to a current running entirely within $\Omega$ will have zero circulation around all loops on
$\partial\Omega$ which bound in $\Omega^\dast$.

We intend to show that the operator \ $\curl\colon\AK(\Omega)\to\K(\Omega)$ is
self-adjoint, and proceed as follows.

Start with a vector field $\vv\in\VF(\Omega)$, and let $\BS(\vv)$ be the 
corresponding magnetic field defined throughout $\reals^3$, $S^3$ or $H^3$.
Let the same symbol denote its restriction to $\Omega$, so that we may 
consider the operator $\BS\colon\VF(\Omega)\to\VF(\Omega)$.
The magnetic field $\BS(\vv)$ is always
divergence-free, but in general is not tangent to the boundary of $\Omega$. 

Now define $\BS'(\vv)$ to be the $L^2$-orthogonal projection of $\BS(\vv)$ into
$\K(\Omega)$. We are only going to apply $\BS'$ to vector fields $\vv$ already
in $\K(\Omega)$, so we regard this \it
modified Biot-Savart operator \rm as a map
$$\BS'\colon\K(\Omega)\to \K(\Omega).$$

We see from the orthogonal decomposition (10.1) that, for any $\vv\in\K(\Omega)$,
we have
$$\BS(\vv)=\BS'(\vv) + \mbox{\rm the gradient component of\ } \BS(\vv).\eqno(11.4)$$

\medskip

\noindent\bf Proposition 11.5\it.\ \ The image of the map $\BS'\colon \K(\Omega)\to\K(\Omega)$
is the subspace $\AK(\Omega)$ of Amp\`erian knots.\rm

\medskip

\noindent\it Proof\rm.\ \ Let $\vv\in\K(\Omega)$, so that $\vv$ is divergence-free and 
tangent to $\partial\Omega$. Then we have\\
$\nabla\times\BS(\vv)=\vv$. This follows
from Maxwell's equation for subdomains of $\reals^3$ by Proposition~1
of Cantarella, DeTurck and Gluck (2001), for subdomains of $S^3$ by 
Proposition~3.1 of
Parsley (2009), and similarly in $H^3$. 

By (11.4), we then also have
$\nabla\times\BS'(\vv)=\vv$.

Now let $\gamma$ be a loop on $\partial\Omega$ which bounds the surface $\Sigma$ in 
$\Omega^\dast$. Then the circulaion of $\BS'(\vv)$ around $\gamma$ equals the flux of 
$\nabla\times\BS'(\vv)=\vv$ through $\Sigma$, according to Amp\`ere's Law. But the flux of 
$\vv$ through $\Sigma$ is zero, since $\vv$ is confined to $\Omega$. Thus $\BS'(\vv)\subset
\AK(\Omega)$, and we have shown that
$$\BS'(\K(\Omega))\subset\AK(\Omega).$$
To see the reverse inclusion, start with $\w\in\AK(\Omega)$, and let
$\vv=\nabla\times\w\in\K(\Omega)$. We claim that $\w=\BS'(\vv)$. 

To show this, first note that $\nabla\times\w=\vv=\nabla\times\BS'(\vv),$
hence $BS'(\vv)-\w$ is curl-free, and therefore lies in $\HK(\Omega)$.

Now we showed above that $\BS'(\vv)\in\AK(\Omega)$, and we have $\w\in\AK(\Omega)$
by hypothesis, so $\BS'(\vv)-\w$ also lies in $\AK(\Omega)$.

Therefore $\BS'(\vv)-\w$ lies in $\AK(\Omega)\cap\HK(\Omega)=\{0\}$, according to
(11.3), and so $\BS'(\vv)=\w$. This shows that
$$\AK(\Omega)\subset\BS'(\K(\Omega)),$$
completing the proof.

\medskip

In the course of the proof, we actually showed a little more.

\medskip

\noindent\bf Corollary 11.6\it.\ \ The maps
$$\BS'\colon\K(\Omega)\to\AK(\Omega)\qquad\mbox{\rm and}\qquad
\curl\colon\AK(\Omega)\to
\K(\Omega)$$
are inverses of one another.\rm

\medskip

\noindent\bf Proposition 11.7\it.\ \ The map $\curl\colon\AK(\Omega)\to \K(\Omega)$
is self-adjoint.\rm

\medskip

\noindent\it Proof\rm.\ \ Self-adjointness of this curl map is equivalent to self-adjointness
of its inverse $\BS'\colon\K(\Omega)\to\AK(\Omega)$, and this in turn is a consequence
of self-adjointness of the Biot-Savart operator $\BS\colon\VF(\Omega)\to\VF(\Omega)$,
which can be seen directly from its 
defining formula as follows. 

\vfill
\eject
Suppose that $\Omega$ is a compact, 
smoothly bounded subdomain of  $\reals^3$, $S^3$ or $H^3$,  and that $\vv,\w\in\VF(\Omega)$.
Then
$$\eqalign{
\ip{\BS(\vv)}{\w}&=\int_\Omega \BS(\vv)(\y)\cdot\w(\y)\,d\y\cr
&=\int_\Omega\left(\int_\Omega P_{\y\x}\vv(\x)\times\nabla_\y\ph_0(\x,\y)\,d\x\right)\cdot \w(\y)\,d\y\cr
&= \int_{\Omega\times\Omega} P_{\y\x}\vv(\x)\times\nabla_\y\ph_0(\x,\y)\cdot\w(\y)\,d\x\,d\y\cr
&= \int_{\Omega\times\Omega} \vv(\x)\times (-\nabla_\x\ph_0(\x,\y))\cdot P_{\x\y}\w(\y)\,d\x\,d\y\cr
&= \int_{\Omega\times\Omega} P_{\x\y}\w(\y)\times\nabla_x\ph_0(\x,\y)\cdot\vv(\x)\,d\x,d\y\cr
&= \int_{\Omega\times\Omega} P_{\y\x}\w(\x)\times\nabla_y\ph_0(\x,\y)\cdot\vv(\y)\,d\x,d\y\cr
&=\vphantom{\int}\ip{\vv}{\BS(\w)},
}$$
where we went from the third line above to the fourth by applying the parallel transport
$P_{\x\y}$ to every term without changing the value of the integrand, from the fourth
to the fifth by interchanging two terms and reversing the sign, from the fifth to the sixth 
by interchanging the variables $\x$ and $\y$, and finally, comparing the sixth line
to the third, moved on to the seventh line. This completes the proof of the proposition.

\medskip

Finally, we come to the desired result.

\medskip

\noindent\bf Theorem 11.8\it.\ \ Let $\Omega$ be a compact, smoothly
bounded subdomain of $\reals^3$, $S^3$ or $H^3$ and let $R=R(\Omega)$ be
the radius of a ball in that space having the same volume as $\Omega$.

Then \ $\curl\colon\AK(\Omega)\to\K(\Omega)$ is a self-adjoint operator,
and for each $\vv\in\AK\Omega)$,
$$|\nabla\times\vv|\ge \frac{|\vv|}{N(R)},$$
where
$$\begin{array}{ll}
\mbox{\it in $\reals^3$ we have}\ &N(R)=R\\ \\
\mbox{\it in $S^3$ we have}\ &N(R)=\displaystyle{\frac{1}{\pi}}(2(1-\cos R)+(\pi-R)\sin R)\\ \\
\mbox{\it in $H^3$ we have}\ &N(R)=\sinh R.
\end{array}$$

\noindent In particular, if $\lambda$ is any curl eigenvalue on $\AK(\Omega)$, then
$\displaystyle{|\lambda|\ge \frac{1}{N(R)}}$.\rm

\medskip

\noindent\it Proof\rm.\ \ We already know from Proposition 11.7 that 
$\curl\colon\AK(\Omega)\to\K(\Omega)$ is self-adjoint.

Let $\vv\in\K(\Omega)$. Then, since $\BS'(\vv)$ is the orthogonal projection
of $\BS(\vv)$ back into the subspace $\K(\Omega)$, we certainly have 
$|\BS'(\vv)|\le|\BS(\vv)|$. 

But $|\BS(\vv)|\le N(R)|\vv|$ by Theorem 9.2, so that same bound applies to
the modified Biot-Savart operator: 
$$|\BS'(\vv)|\le N(R)|\vv|.$$ 
Now suppose that $\vv\in\AK(\Omega)$. Then
$$\nabla\times\vv\in\K(\Omega)\qquad\mbox{\rm and}\qquad\BS'(\nabla\times\vv)=\vv.$$
Hence
$$|\vv|=|\BS'(\nabla\times\vv)|\le N(R)|\nabla\times\vv|,$$
and therefore
$$|\nabla\times\vv|\ge\frac{1}{N(R)}|\vv|,$$
completing the proof of Theorem 11.8.

\bigskip

\noindent\bf References\rm

\addtolength{\baselineskip}{-2pt}

\begin{description}
\item[1820] Jean-Baptiste Biot and Felix Savart, {\it Note sur le magnetisme de la pile de Volta},
Annales de chimie et de physique, 2nd ser.,  \bf 15\rm, 222--223.

\item[1824] Jean-Baptiste Biot, {\it Precise Elementaire de Physique Experimentale},
3rd edition, Vol II, Chez Deterville (Paris).

\item[1833] Carl Friedrich Gauss, {\it Integral formula for linking number}, in
{\it Zur mathematischen theorie der electrodynamische wirkungen},
Collected Works, Vol 5, K\"oniglichen Gesellschaft des Wissenschaften,
Gottingen, 2nd edition, page 605.

\item[1891] James Clerk Maxwell, {\it A Treatise on Electricity and Magnetism},
Reprinted by the Clarendon Press, Oxford (1998), two volumes.

\item[1958] L. Woltjer, {\it A theorem on force-free magnetic fields},
Proc.\ Nat.\ Acad.\ Sci.\ USA \bf 44\rm, 489--491.

\item[1959] Georges C{\accent21 a}lug{\accent21 a}reanu, {\it L'integral de Gauss et l'analyse
des n\oe uds tridimensionnels},
Rev.\  Math.\  Pures Appl.\  \bf 4\rm, 5--20.

\item[1960] Georges de Rham, {\it Vari\'et\'es Différentiables}, 2nd edition, Hermann, Paris.

\item[1961] Georges C{\accent21 a}lug{\accent21 a}reanu, {\it Sur les classes d'isotopie des n\oe
uds tridimensionnels et
leurs invariants}, Czechoslovak Math.\  J.\  \bf 11\rm(86), 588--625.

\item[1961] Georges C{\accent21 a}lug{\accent21 a}reanu, {\it Sur les enlacements tridimensionnels des courbes ferm\'ees},
Comm.\ Acad.\ R\. P.\  Romine \bf 11\rm, 829--832.

\item[1968a]  William Pohl, {\it The self-linking number of a closed space curve}, Journal of
Mathematics and Mechanics, \bf 17\rm(10), 975--985.

\item[1968b] William Pohl, {\it Some integral formulas for space curves and their generalizations},
American Journal of Mathematics \bf 90\rm, 1321--1345.

\item[1969] H. K. Moffatt, {\it The degree of knottedness of tangled vortex lines},
J.\  Fluid Mech.\  \bf 35\rm, 117--129 and 159, 359--378.

\item[1969]  James White, {\it Self-linking and the Gauss integral in higher dimensions},
American Journal of Mathematics \bf 91\rm, 693--728.
\item[1971] F.\  Brock Fuller, {\it The writhing number of a space curve},
Proc.\  Nat.\  Acad.\  Sci.\  USA  \bf 68\rm(4), 815--819.

\item[1974] V.\  I.\  Arnold, {\it The asymptotic Hopf invariant and its applications},
English translation in Selecta Math.\  Sov., \bf 5\rm(4), (1986) 327--342;
original in Russian, Erevan (1974).

\item[1981] David Griffiths, {\it Introduction to Electrodynamics}, Prentice Hall, New Jersey.
second edition 1989, third edition 1999.

\item[1990] Robert Zimmer, {\it Essential Results of Functional Analysis\rm}, University of 
Chicago Press, Chicago and London.  

\item[1992] H.\  K.\  Moffatt and R.\  L.\  Ricca, {\it Helicity and the C{\accent21 a}lug{\accent21 a}reanu invariant},
Proc.\  Royal Soc.\  Lond.\  A, \bf 439\rm, 411--429.

\item[1992] R.\  L.\  Ricca and H.\  K.\  Moffatt, {\it The helicity of a knotted vortex filament},
Topological Aspects of Dynamics of Fluids and Plasmas (H.\  K.\  Moffatt,
et al, eds.) Kluwer Academic Publishers (Dordrecht, Boston), 225--236.

\item[1995] Gerald Folland, {\it Introduction to Partial Differential Equations}, 
Princeton University Press, Princeton.

\item[1998] Moritz Epple, {\it Orbits of asteroids, a braid, and the first link invariant},
The Mathematical Intelligencer, \bf 20\rm(1)

\item[1999] Jason Cantarella, {\it Topological structure of stable plasma flows},
Ph.D.\  thesis, University of Pennsylvania.

\item[1999] Jason Cantarella, Dennis DeTurck, Herman Gluck and Mikhail Teytel,
{\it Influence of geometry and topology on helicity}, Magnetic Helicity in Space  
and Laboratory Plasmas, ed.\  by M.\  Brown, R.\  Canfield and A.\  Pevtsov,
Geophysical Monograph, \bf 111\rm, American Geophysical Union, 17--24.

\item[2000a] Jason Cantarella, Dennis DeTurck, Herman Gluck and Mikhail
Teytel,
{\it Isoperimetric problems for the helicity of vector fields and the Biot-Savart
and curl operators}, Journal of Mathematical Physics  \bf 41\rm(8), 5615--5641.

\item[2000b] Jason Cantarella, Dennis DeTurck and Herman Gluck, {\it
Upper bounds for the
writhing of knots and the helicity of vector fields}, Proc.\  of Conference in Honor
of the 70th Birthday of Joan Birman, ed. by J.\  Gilman, X.-S.\  Lin and W.\  Menasco,    
International Press, AMS/IP Series on Advanced Mathematics.

\item[2000c] Jason Cantarella, Dennis DeTurck, Herman Gluck and Misha
Teytel,
{\it Eigenvalues and eigenfields of the Biot-Savart and curl operators on
spherically symmetric domains}, Physics of Plasmas \bf 7\rm(7), 2766--2775.

\item[2001] Jason Cantarella, Dennis DeTurck and Herman Gluck, {\it The Biot-Savart
operator for application to knot theory, fluid dynamics and plasma physics},
Journal of Mathematical Physics \bf 42\rm(2), 876--905.

\item[2002] Jason Cantarella, Dennis DeTurck and Herman Gluck, {\it Vector calculus and the
topology of domains in 3-space}, American Mathematical Monthly \bf 109\rm(5), 409--442.

\item[2004] Jason Parsley, {\it The Biot-Savart operator and electrodynamics on bounded subdomains
of the 3-sphere}, Ph.D.\ thesis, University of Pennsylvania.

\item[2004] Dennis DeTurck and Herman Gluck, {\it The 
Gauss Linking Integral on the $3$-sphere
and in hyperbolic $3$-space}, arXiv:math/0406276v1 [math.GT].

\item[2008] Dennis DeTurck and Herman Gluck, {\it Electrodynamics and the 
Gauss Linking Integral on the $3$-sphere
and in hyperbolic $3$-space}, Journal of Mathematical Physics \bf 49\rm(2), 023504.

\item[2008] Greg Kuperberg, {\it From the Mahler conjecture to Gauss linking forms}, 
Geom.\ Funct.\ Anal.\ \bf 18\rm(3), 870--892.

\item[2009] Jason Parsley, {\it The Biot-Savart operator and electrodynamics on subdomains
of the $3$-sphere}, arXiv:0904.3524v1 [math.DG].

\end{description}

\noindent University of Pennsylvania\\
Philadelphia, PA  19104

\noindent
\it deturck@math.upenn.edu\rm\\
\it gluck@math.upenn.edu\rm\\

\end{document}        

\end{document}